\documentclass[11pt,reqno]{amsart}
\usepackage{amsmath,amssymb}
\usepackage{array}
\usepackage[svgnames]{xcolor}
\usepackage{tikz}
\usetikzlibrary{shapes.geometric}
\usepackage{amsfonts} 
\usepackage{mathtools}
\usepackage{enumerate}
\usepackage{caption}
\usepackage{hyperref}
\usepackage{doi}
\hypersetup{colorlinks,linkcolor=blue ,citecolor=blue, urlcolor=blue}

\numberwithin{equation}{section}  

\DeclareMathOperator{\aff}{aff}
\DeclareMathOperator{\Span}{span}

\DeclareMathOperator{\conv}{conv}
\DeclareMathOperator{\id}{id}

\newcommand{\filter}[1]{\langle #1 \rangle_\uparrow}
\newcommand{\ideal}[1]{\langle #1 \rangle_\downarrow}

\def\RR{\mathbb{R}}
\def\CP{\mathcal{C}(P)}
\def\OP{\mathcal{O}(P)}
\def\Fc{{\mathcal F}}
\def\Cc{{\mathcal C}}
\def\Pc{{\mathcal P}}
\def\Ac{{\mathcal A}}
\def\Oc{{\mathcal O}}

\newtheorem{theorem}{Theorem}[section]
\newtheorem{lemma}[theorem]{Lemma}
\newtheorem{corollary}[theorem]{Corollary}
\newtheorem{proposition}[theorem]{Proposition}
\newtheorem{remark}[theorem]{Remark}
\newtheorem{definition}[theorem]{Definition}
\newtheorem{conjecture}[theorem]{Conjecture}
\newtheorem{setup}[theorem]{Setup}

\begin{document}

\title{Two-dimensional faces of order and chain polytopes}
\author{Ragnar Freij-Hollanti, Teemu Lundström, Aki Mori}

\address{Ragnar Freij-Hollanti, Department of Mathematics and Systems Analysis, Aalto University, Espoo, Finland}
\email{ragnar.freij@aalto.fi}
\address{Teemu Lundström, Department of Mathematics and Systems Analysis, Aalto University, Espoo, Finland}
\email{teemu.lundstrom@gmail.com}
\address{Aki Mori, Center for Physics and Mathematics, Institute for Liberal Arts and Sciences, Osaka Electro-Communication University, Neyagawa, Osaka 572-8530, Japan}
\email{a-mori@osakac.ac.jp}

\subjclass[2020]{52B05, 06A07}
\keywords{order polytope, chain polytope, partially ordered set, $f$-vector}

\begin{abstract}
We give an explicit combinatorial description of the $2$-dimensional faces of both the order polytope $\Oc(P)$ and the chain polytope $\Cc(P)$ of a partially ordered set $P$. Using these descriptions, we show that for any $P$, $\Cc(P)$ has equally many square faces, and at least as many triangular faces, as $\Oc(P)$ does. Moreover, the inequality is shown to be strict except when $\Oc(P)$ and $\Cc(P)$ are unimodularly equivalent. This proves the case $i=2$ of a conjecture by Hibi and Li.
\end{abstract}

\maketitle

\tableofcontents

\section{Introduction}
A {\it 0/1-polytope} $\Pc \subseteq \RR^d$ is a convex polytope whose vertices all have coordinates in $\{0,1\}$.
Stanley \cite{Sta} introduced two types of $0/1$-polytopes associated with a finite partially ordered set $P$ (``poset'' for short). These are known as the order polytope $\Oc(P)$ and the chain polytope $\Cc(P)$.
They form important classes of $d$-dimensional polytopes that have been widely studied in combinatorics and commutative algebra. The order polytope $\Oc(P)$ had previously been introduced by Geissinger in~\cite{Ge}.
In this work, we study the face structure and the number of faces of $\Oc(P)$ and $\Cc(P)$.
Let $f_i = f_i(\Pc)$ denote the number of $i$-faces of $\Pc$. The integer sequence $f(\Pc) = (f_0, f_1, \dots, f_{d-1})$ is called the {\it $f$-vector} of $\Pc$.
Hibi and Li proposed the following conjecture with respect to the $f$-vectors of $\Oc(P)$ and $\Cc(P)$:

\begin{conjecture}[{\cite[Conjecture 2.4]{HL2}}]\label{HL-conj}
Let $P$ be a finite poset with $|P| = d > 1$. 
Then 
\begin{enumerate}[(a)]
\item $f_i(\Oc(P)) \leq f_i(\Cc(P))$ for all $1 \leq i \leq d-1$, and
\item if $f_i(\Oc(P)) = f_i(\Cc(P))$ for some $2 \leq i \leq d-1$, then $\Oc(P)$ and $\Cc(P)$ are unimodularly equivalent.    
\end{enumerate}
\end{conjecture}

\begin{figure}
    \centering
    \begin{tikzpicture}
    [scale=0.6,
    dot/.style={draw, circle, minimum size=2mm, inner sep=0pt, fill=white},
    ]
\node(a)[dot] at (0,0){};
\node(b)[dot] at (2,0){};
\node(c)[dot] at (1,1.5){};
\node(d)[dot] at (0,3){};
\node(e)[dot] at (2,3){};

\draw (a) -- (c);
\draw (b) -- (c);
\draw (d) -- (c);
\draw (e) -- (c);
\end{tikzpicture}
    \caption{The smallest poset $\mathbf{X}$ for which $\Oc(\mathbf{X})$ and $\Cc(\mathbf{X})$ are not unimodularly equivalent.}\label{fig:X}
\end{figure}
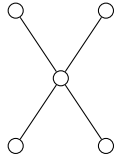

It is shown in \cite{Sta} that $f_0(\Oc(P)) = f_0(\Cc(P))$.
In \cite{HL2}, it was proved that $f_{d-1}(\Oc(P)) \le f_{d-1}(\Cc(P))$ with equality holding if and only if $\Oc(P)$ and $\Cc(P)$ are unimodularly equivalent, and that this happens if and only if $P$ is {\em $\mathbf{X}$-free}, meaning that it does not contain the poset $\mathbf{X}$ depicted in Figure~\ref{fig:X} as a subposet.
Relatedly, \cite{BCMOS} focuses on the difference $f_{d-1}(\Cc(P)) - f_{d-1}(\Oc(P))$ and classifies the posets for which the difference is exactly one.
Moreover, it was shown in \cite{HLSS} that $f_1(\Oc(P)) = f_1(\Cc(P))$ for any poset $P$.  
Two-dimensional faces of $\mathcal{O}(P)$ and $\mathcal{C}(P)$ have been studied in \cite{J}.

In \cite{FL}, it was proved that if $P$ is a poset built inductively by taking disjoint unions and ordinal sums of posets, starting from $\mathbf{X}$-free posets, then $f_i(\Oc(P)) \le f_i(\Cc(P))$ holds for all $1 \le i \le d - 1$.
In \cite{M2}, it was shown that if $P$ is a {\em maximal ranked} poset, then the number of triangular 2-faces of $\Oc(P)$ is less than or equal to that of $\Cc(P)$, with equality holding if and only if $\Oc(P)$ and $\Cc(P)$ are unimodularly equivalent. In this paper, we generalize this result to the class of all posets.

The study of {\em marked} order and chain polytopes, which generalize the concepts of $\Oc(P)$ and $\Cc(P)$, was initiated in \cite{ABS}, and has since been actively studied, motivated by connections to representation theory of Lie algebras.
A generalization of Conjecture~\ref{HL-conj}(a) to the setting of marked chain-order polytopes was proposed in \cite[Conjecture~5.3]{FF}, and it was confirmed in \cite{AFJ} for the case where $P$ is a (unmarked) maximal ranked poset.  

In the present paper, we focus on the case of $2$-dimensional faces in Conjecture~\ref{HL-conj}.  
Throughout this paper, we refer to a $2$-dimensional face of a convex polytope as a {\em square face} if it has four vertices, and as a {\em triangular face} (or a triangle) if it has three vertices. We will denote by $f_2^\Box(\Pc)$ the number of square faces and by $f_2^\Delta(\Pc)$ the number of triangle faces respectively, of a polytope $\Pc$. It is an easy observation (Corollary~\ref{cor:triangle_or_square}) that every 2-dimensional face of $\Oc(P)$ and $\Cc(P)$ is a triangle or a square, so $$f_2(\Oc(P))=f_2^\Delta(\Oc(P))+f_2^\Box(\Oc(P))\mbox{ and }f_2(\Cc(P))=f_2^\Delta(\Cc(P))+f_2^\Box(\Cc(P)).$$

The main results of this paper are summarized as follows.
\begin{itemize}
  \item New combinatorial parametrizations of the triangular faces of $\Oc(P)$ (Proposition~\ref{OT_Char}) and $\Cc(P)$ (Proposition~\ref{CT_Char}).
  \item $f_2^\Delta(\Oc(P))\leq f_2^\Delta(\Cc(P))$ for any poset $P$, with equality holding if and only if $P$ is $\mathbf{X}$-free (Theorem~\ref{thm:triangle_ineq}).
  \item Explicit geometric descriptions of square faces of $\Oc(P)$ (Theorem~\ref{SQorder}) and $\Cc(P)$ (Theorem~\ref{thm:square_faces_of_C(P)}).
  \item $f_2^\Box(\Oc(P))= f_2^\Box(\Cc(P))$ for any poset $P$ (Theorem~\ref{thm:square_eq}).
\end{itemize}
As a consequence, this paper proves the case $i=2$ of Conjecture~\ref{HL-conj}. As a corollary, we also see that $\mathcal{O}(P)$ and $\mathcal{C}(P)$ are unimodularly equivalent if and only if $f_2(\mathcal{O}(P)) = f_2(\mathcal{C}(P))$.
In fact, unimodularity holds if and only if $f_2^\Delta(\mathcal{O}(P)) = f_2^\Delta(\mathcal{C}(P))$.

The present paper is organized as follows.  
In Section~\ref{basics}, we review basics on posets, along with known explicit descriptions of the faces of $\Oc(P)$ and $\Cc(P)$.  
In Section~\ref{Preliminary}, we present preliminary results on polytopes and posets.
In particular, in Section~\ref{geometric_lemmas} we show that every $2$-dimensional face of a $0/1$-polytope is either a triangle or a square. 
In Section~\ref{subsection:recharacterize_triangles}, we give a new combinatorial description for the triangle faces in both $\mathcal{O}(P)$ and $\mathcal{C}(P)$.
In Section~\ref{subsection:triangle_ineq} we develop recursive (under quotients of subposets) combinatorial formulas for the number of triangles in $\Oc(P)$ and $\Cc(P)$. 
These formulas are then used to inductively prove that the number of triangular faces is always smaller or equal in $\Oc(P)$, than in $\Cc(P)$.
In Section~\ref{subsection:square_of_order}, we provide an explicit description of square faces of $\Oc(P)$ in terms of filters.  
In Section~\ref{subsection:square_chain}, we give an explicit description of square faces of $\Cc(P)$ in terms of antichains. 
In this section, we also study the poset consisting of the set of antichains of $P$, ordered by inclusion, and show that this poset reflects certain types of edges in $\Cc(P)$.  
In Section~\ref{subsection:Equality}, we establish the equality of the number of square faces of $\Oc(P)$ and $\Cc(P)$.

\section{The two polytopes and other definitions}\label{basics}

\subsection{Notation for partially ordered sets}
In this section, we collect the basic terminology that we will use for posets. The experienced reader may want to skip this section at a first read and refer back as needed. 

Let $P$ be a finite poset equipped with a partial order $\leq$. 
We say that $x$ and $y$ are {\em comparable} in $P$, denoted by $x \perp y$, if $x \le y$ or $y \le x$.
Otherwise they are {\em incomparable}, denoted by $x \parallel y$.
Let $A$ and $B$ be subsets of $P$. We say $A$ and $B$ are {\em incomparable}, denoted $A\parallel B$, if $a\parallel b$ for all $a\in A$ and $b\in B$. 
Note that $A \parallel B$ implies $A \cap B = \emptyset$.
We denote by $P^\mathrm{op}$ the opposite poset of $P$, which is the poset whose ground set is identical with that of $P$, and with $x\le y$ in
$P^\mathrm{op}$ if and only if $y\le x$ in $P$.

A non-empty poset $P$ is said to be {\it connected} if for every $x, y \in P$, there exists a sequence 
$$
x= p_1 \perp p_2 \perp  \cdots \perp p_s = y
$$
where $p_1, \ldots, p_s \in P$. 
Equivalently, $P$ is connected if it is not the disjoint union $P=Q\sqcup R$ of two non-empty subposets with $Q\parallel R$.

A {\em filter} (or an {\em upset}) of $P$ is a subset $F \subseteq P$ such that if $x \in F$ and $y \ge x$, then $y \in F$.
Dually, an {\em order ideal} (or a {\em downset}) of $P$ is a subset $I \subseteq P$ such that if $x \in I$ and $y \le x$, then $y \in I$.
An {\em antichain} of $P$ is a subset $A \subseteq P$ such that $x$ and $y$ belonging to $A$ with $x \neq y$ are incomparable.
Note that the empty set $\emptyset$ is considered both a filter and an antichain of~$P$.

For a subset $S\subseteq P$ we get two antichains $\min (S)$ and $\max (S)$ consisting of the minimal and maximal elements of $S$ respectively. We also get a filter 
$$
\langle S \rangle_\uparrow:= \{ y \in P : y \geq x \text{ for some } x \in S\},
$$ and an order ideal $$
\langle S \rangle_\downarrow:= \{ y \in P : y \leq x \text{ for some } x \in S\}.
$$
There is a bijection between filters $F$ and antichains $A$ of $P$, given by
\begin{gather*}
F = \filter{A}\quad\text{and}\quad
A = \min(F).
\end{gather*}
Similarly, there is a bijection between order ideals and antichains.
For a given $x \in P$, the set 
$\filter{x} \coloneqq \filter{\{ x \}}= \{y \in P : y \ge x\}$ is called the {\it principal filter} generated by
$x$.
Similarly, the set $\ideal{x} \coloneqq \ideal{\{ x \}} = \{y \in P : y \le x\}$ is called {\it principal order ideal} generated by $x$.

A {\em chain} in $P$ is a set of pairwise comparable elements in $P$, and the {\em height} of $P$ is the maximum size of a chain in $P$.
We say that $x$ {\it covers} $y$ if $x < y$ and $x < z < y$ for no $z \in P$. 
A subset $S \subseteq P$ is called \emph{order-convex} if $x,y \in S$ and $x \le z \le y$ implies $z \in S$.

In Section~\ref{Triangles} we will use contractions of posets.
Contractions of posets often need to be defined with care, see {\em e.g.}~\cite{Will}, but in our case we will only need to do contractions by either ideals or filters, where the definitions are more straightforward.

\begin{definition}
    Let $P$ be a poset and let $J$ be either a filter or an order ideal.
    Partition $P$ as $\{ \{ x \} \colon x \in P \setminus J \} \cup \{ J \}$ and let $[p]$ denote the part containing $p \in P$.
    We define a relation $\le_J$ on these parts by setting $[x] \le_J [y]$ if $x' \le y'$ for some $x' \in [x]$ and $y' \in [y]$.
    By $P / J$ we mean the poset with ground set $\{ [x] \colon x \in P \}$ and order relation $\le_J$.
\end{definition}
The above definition relies on the fact that $\le_J$ is indeed a partial order on the set $\{ [x] \colon x \in P \}$. 
Checking that this holds when $J$ is a filter or an order ideal is straightforward and left to the reader.

If $J$ is clear from the context, we shall write $\le$ instead of $\le_J$.
For any subset $S \subseteq P$ we write $S / J \coloneqq \{ [x] \colon x \in S \} \subseteq P/J$.
If $J \neq \emptyset$, we will denote the set $[q]$, $q \in J$, by $v_J$.

\begin{remark}\label{rem:identify}
    Let $P$ be a poset and $J$ a filter or an order ideal.
    Let $S \subseteq P$ be any subset.
    \begin{enumerate}[(1)]
        \item If $S \cap J = \emptyset$ then $S$ can be identified with $S / J$ via $x \leftrightarrow \{ x \}$. 
        We shall identify $S$ with $S/J$ when $S \cap J = \emptyset$ without further comment.
        \item $S \cap J \neq \emptyset$ if and only if $v_J \in S/J$.
        \item If $J \neq \emptyset$ is a filter (ideal) then $v_J$ is maximal (minimal) in $P/J$.
        \item We may think of $P/J$ as the poset with ground set $(P \setminus J) \cup \{ v_J \}$  where $v_J$ is a new minimal or maximal element.
    \end{enumerate}
\end{remark}

\subsection{Order and chain polytopes}
Stanley \cite{Sta} defined two types of lattice polytopes arising from a poset.
The {\em order polytope} of $P$ is defined by the following:
\begin{equation*}
  \Oc(P) 
  =
  \left\{ x \in \mathbb{R}^P \;\middle|\;
  \begin{gathered} 
    0 \leq x_p \leq 1 \mbox{ for all } p \in P, \\
    x_p \leq x_q \mbox{ if } p \leq q \mbox{ in }P
  \end{gathered}
  \right\}
  \subseteq \mathbb{R}^P.
\end{equation*}
The {\em chain polytope} of $P$ is defined by the following:
\begin{equation*}
  \Cc(P) 
  =
  \left\{ x \in \mathbb{R}^P \;\middle|\;
  \begin{gathered} 
    x_p \geq 0 \mbox{ for all } p \in P,\\
    x_{p_1}+\cdots +x_{p_k} \leq 1 \mbox{ if } p_1 < \cdots < p_k \mbox{ in } P
  \end{gathered}
  \right\}
  \subseteq \mathbb{R}^P.
\end{equation*}

We briefly review the properties of these two polytopes as established by \cite{Sta}.
One has 
$$
{\rm dim}(\Oc(P)) = {\rm dim}(\Cc(P))=|P|.
$$
To each subset $W \subseteq P$, we associate $\chi_W \in \mathbb{R}^P$ which is the characteristic vector (function) of $W$. 
In particular $\chi_\emptyset$ is the origin of $\mathbb{R}^P$.
The vertex sets of $\Oc(P)$ and $\Cc(P)$ are given as follows, respectively:
\begin{gather*}
{\rm vert}({\Oc}(P))
=
\{\chi_F : F \mbox{ is a filter of } P\},\\
{\rm vert}({\Cc}(P))
=
\{\chi_A : A \mbox{ is an antichain of } P\}.
\end{gather*}
By the bijections between antichains and filters, it follows that the number of vertices of $\Oc(P)$ is equal to that of $\Cc(P)$. 

The facet-defining inequalities of $\mathcal{O}(P)$ are
\begin{itemize}
\item 
$x_p \ge 0$, for all $p \in P$ minimal,
\item
$x_q \le 1$, for all $q \in P$ maximal, and
\item
$x_p \le x_q$, for all pairs $(p,q)$ where $q$ covers $p$.
\end{itemize}
\noindent
For $\mathcal{C}(P)$, the facet defining inequalities are

\begin{itemize}
\item 
$x_p \ge 0$, for all $p  \in P$, and
\item
$x_{p_1} + \dots + x_{p_k} \le 1$, for all maximal chains $p_1 < \cdots < p_k$ of $P$.
\end{itemize}

The descriptions of edges of $\Oc(P)$ and $\Cc(P)$ are obtained as follows.
Here $A \Delta B$ denotes symmetric difference of $A$ and $B$, that is $A \Delta B = (A \setminus B)\cup(B \setminus A)$. 

\begin{proposition}[{\cite[Lemma 4, Lemma 5]{HL1}}]\label{edge}
Let $P$ be a finite poset.

\begin{enumerate}
\item[(a)]
Let $F$ and $G$ be filters of $P$ with $F \neq G$. 
Then $\conv(\chi_F,\chi_G)$ is an edge of $\Oc(P)$ if and only if $F \subseteq G$ and $G \setminus F$ is connected in $P$.
\item[(b)]
 Let $A$ and $B$ be antichains of $P$ with $A \neq B$.
Then $\conv(\chi_A,\chi_B)$ is an edge of $\Cc(P)$ if and only if $A \Delta B$ is connected in $P$.
\end{enumerate}
\end{proposition}

The descriptions of triangular faces of $\Oc(P)$ and $\Cc(P)$ are obtained as follows.

\begin{proposition}[{\cite[Lemma 2.2]{M2}}]\label{triangle face}
Let $F$, $G$, and $H$ be pairwise distinct filters of $P$ and let $A$, $B$, and $C$ be pairwise distinct antichains of $P$.
\begin{enumerate}
\item[(a)]
$\conv(\chi_F, \chi_G, \chi_H)$ is a 2-face of $\Oc(P)$ if and only if $F \subseteq G \subseteq H$ and $G \setminus F$, $H \setminus G$, $H \setminus F$ are connected in $P$.
\item[(b)]
$\conv(\chi_A, \chi_B, \chi_C)$ is a 2-face of $\Cc(P)$ if and only if $A \Delta B$, $B \Delta C$ and $C \Delta A$ are connected in $P$.
\end{enumerate}
\end{proposition}

Proposition \ref{triangle face} implies that the triangles in the $1$-skeleton of $\Oc(P)$ or $\Cc(P)$ correspond to the triangular faces of each polytope.

\section{Preliminaries}\label{Preliminary}
In this section, we give some preliminary results that will be helpful throughout the paper.
In the first subsection we focus on geometric results and in the second we prove various lemmas related to posets. The experienced reader can skip this section at the first read, and refer back to it as necessary.

\subsection{Geometric lemmas}\label{geometric_lemmas}
We start by showing that there are only two possible combinatorial types of $2$-dimensional faces of order polytopes and chain polytopes: a triangle and a square.  
This is obtained from a more general result which shows that any $0/1$-polytope is isomorphic to a full-dimensional $0/1$-polytope.
In fact, we will show an even more general result that any polyhedron is isomorphic to a full-dimensional polyhedron via coordinate projections.

For all $i \in [n] = \{ 1,\dots,n \}$ we let 
\begin{equation*}
    \pi_i \colon \mathbb{R}^n \to \mathbb{R}^{n-1} , \ \pi_i(x_1, \dots, x_n) = (x_1,\dots,x_{i-1},x_{i+1},\dots,x_n).
\end{equation*}
In what follows, $\Span$ refers to linear span, $\aff$ refers to affine span and $\conv$ refers to convex hull.

\begin{lemma}\label{lemma:intersecting_coordinate_lines_implies_full_dimensionality}
    If $S \subseteq \mathbb{R}^n$ is any non-empty set such that for all $i \in [n]$ there exists distinct $x,y \in S$ with $x-y \in \Span(e_i)$, then $S$ is full-dimensional.
\end{lemma}

\begin{proof}
    Fix $x_0 \in S$.
    We have 
    \begin{equation*}
        \aff(S) = x_0 + \Span(x-y \colon x,y \in S).
    \end{equation*}
    By assumption, for all $i \in [n]$ there exists distinct $x,y \in S$ such that $x-y = \lambda e_i$ for some $\lambda \neq 0$.
    Hence $\Span(x-y \colon x,y \in S)$ contains all standard basis vectors and hence $S$ is full-dimensional.
\end{proof}

\begin{corollary}\label{cor:injX}
    If $X \subseteq \mathbb{R}^n$ is a non-empty and non-full-dimensional set, then there exists $i \in [n]$ such that $\pi_i$ is injective on $\aff(X)$.
\end{corollary}

\begin{proof}
    Since $\aff(X)$ is not full-dimensional, by Lemma \ref{lemma:intersecting_coordinate_lines_implies_full_dimensionality} there exists $i \in [n]$ such that for all distinct $x,y \in \aff(X)$ we have  $x-y \not\in \Span(e_i)$.
    Since $\Span(e_i) = \ker(\pi_i)$ we have that $\pi_i$ is injective on $\aff(X)$.
\end{proof}

Note that since $\pi_i$ is injective on $\aff(X)$ it follows that $\dim(X) = \dim(\pi_i(X))$.
In general, if an affine map $f$ is injective on $\aff(X)$ then $\dim(X) = \dim(f(X))$.

\begin{proposition}\label{prop:projection_injecive}
    For any non-empty set $X \subseteq \mathbb{R}^n$ with $\dim(X) = d$ there exists a coordinate projection $f \colon \mathbb{R}^n \to \mathbb{R}^d$ such that $f$ is injective on $\aff(X)$.
\end{proposition}

\begin{proof}
    We prove the proposition with induction on $n-d$.
    In the base case $n=d$ in which case we can choose $f = \id$.
    Suppose $d < n$.
    By Corollary~\ref{cor:injX} we find $i \in [n]$ such that $\pi_i \colon \mathbb{R}^n \to \mathbb{R}^{n-1}$ is injective on $\aff(X)$.
    Since $\dim(\pi_i(X)) = d \le n-1$, by induction assumption we find a coordinate projection $f \colon \mathbb{R}^{n-1} \to \mathbb{R}^d$ such that $f$ is injective on $\aff(\pi_i(X)) = \pi_i(\aff(X))$.
    Now $f \circ \pi_i \colon \mathbb{R}^n \to \mathbb{R}^d$ is a coordinate projection that is injective on $\aff(X)$.
\end{proof}

\begin{proposition}\label{prop:isom_to_full_dim}
    Any non-empty polyhedron is affinely isomorphic to a full-dimensional polyhedron via coordinate projection.
    In particular, any $0/1$-polytope is affinely isomorphic to a full-dimensional $0/1$-polytope.
\end{proposition}

\begin{proof}
    Let $P \subseteq \mathbb{R}^n$ be any non-empty polyhedron where $\dim(P) = d$.
    By Proposition~\ref{prop:projection_injecive} we find a coordinate projection $f \colon \mathbb{R}^n \to \mathbb{R}^d$ such that $f$ is injective on $\aff(P)$.
    Thus $P \cong f(P)$ where $\dim(f(P)) = \dim(P) = d$.
    
    The latter claim follows from the fact that if $P = \conv(v_1,\dots,v_k)$ is a 0/1-polytope then $f(P) = \conv(f(v_1),\dots,f(v_k))$ is also a 0/1-polytope.
\end{proof}

Since the only two-dimensional 0/1-polytopes in $\mathbb{R}^2$ are triangles and squares, Proposition~\ref{prop:isom_to_full_dim} immediately implies the following.

\begin{corollary}\label{cor:triangle_or_square}
    A 2-face of $\mathcal{O}(P)$ or $\mathcal{C}(P)$ is either a triangle or a square.
\end{corollary}

By Corollary~\ref{cor:triangle_or_square}, we can write $f_2(\Pc)=f_2^\Delta(\Pc)+f_2^\Box(\Pc)$ when $\Pc$ is an order or a chain polytope, with notation as in the introduction. More generally, by Proposition~\ref{prop:isom_to_full_dim}, this holds for any $0/1$-polytope.

In Sections~\ref{subsection:square_of_order} and \ref{subsection:square_chain} we will study the square faces of $\mathcal{O}(P)$ and $\mathcal{C}(P)$ which we will characterize in terms of filters and antichains of $P$.
We will approach this by looking at how four filters or four antichains corresponding to a square face can be contained in each other.
The following results will help us rule out some possible containments between such filters and antichains.

\begin{lemma}\label{Lemma:s(i)}
    Let $A,B,C,D \subseteq [n]$ be distinct subsets.
    If some $i \in [n]$ is contained in exactly 1 or 3 of the set $A,B,C,D$, then 
    \begin{equation*}
        \dim(\conv(\chi_A,\chi_B,\chi_C,\chi_D)) = 3.
    \end{equation*}
\end{lemma}

\begin{proof}
    Say $i \in (A \cap B \cap C) \setminus D$.
    Now the 0/1-vectors $\chi_A,\chi_B,\chi_C$ lie on the hyperplane $\{ x \in \mathbb{R}^n \mid x_i = 1 \}$ while $\chi_D$ does not.
    Thus $\dim(\conv(\chi_A,\chi_B,\chi_C,\chi_D)) = 3$.
    Then suppose $i \in D \setminus (A \cup B \cup C)$.
    Now the vectors $\chi_A,\chi_B,\chi_C$ are on the hyperplane $\{ x \in \mathbb{R}^n \mid x_i = 0 \}$ while $\chi_D$ is not.
    Hence in this case also $\dim(\conv(\chi_A,\chi_B,\chi_C,\chi_D)) = 3$.
\end{proof}

\begin{corollary}\label{cor:onlyone}
    If $A,B,C,D \subseteq [n]$ are distinct sets such that $$
    \dim(\conv(\chi_A,\chi_B,\chi_C,\chi_D)) = 2
    $$
    then there can be no element in $[n]$ that belongs to exactly 1 or  3 of the sets $A,B,C,D$.
\end{corollary}

\begin{lemma}\label{Lemma:princip3}
    Let $A,B,C,D \subseteq [n]$ be distinct subsets.
    If in the poset $(\{ A,B,C,D \},\subseteq)$ there is a principal filter or a principal order ideal of size 3 then 
    \begin{equation*}
        \dim(\conv(\chi_A,\chi_B,\chi_C,\chi_D)) = 3.
    \end{equation*}
\end{lemma}

\begin{proof}
    Suppose $A$ generates a filter, say $\{ A,B,C \}$.
    Here $A \subseteq B$, $A \subseteq C$ but $A \not\subseteq D$.
    Thus there exists $i \in A \setminus D$.
    Hence $i \in (A \cap B \cap C) \setminus D$.
    Thus $\dim(\conv(\chi_A,\chi_B,\chi_C,\chi_D)) = 3$ by Lemma \ref{Lemma:s(i)}.
    Then suppose $A$ generates an order ideal, say $\{ A,B,C \}$.
    Here $B \subseteq A$, $C \subseteq A$ and $D \not\subseteq A$.
    Hence there exists $i \in D \setminus A$.
    Thus $i \in D \setminus (A \cup B \cup C)$.
    Again, $\dim(\conv(\chi_A,\chi_B,\chi_C,\chi_D)) = 3$ by Lemma \ref{Lemma:s(i)}.
\end{proof}

The last result of this subsection shows how the sets corresponding to the diagonals of a 0/1-square relate to each other.

\begin{corollary}\label{cor:diagonals}
    If $A,B,C,D \subseteq [n]$ are distinct such that $S = \conv(\chi_A,\chi_B,\chi_C,\chi_D)$ is a square where the diagonals are $\conv(\chi_A,\chi_D)$ and $\conv(\chi_B,\chi_C)$, then $\chi_A + \chi_D = \chi_B + \chi_C$.
    In particular,  $A \cap D = B \cap C$ and $A \cup D = B \cup C$.
\end{corollary}

\begin{proof}
    If $n = 2$ then the claim is clear.
    If $n > 2$ then by Corollary \ref{prop:projection_injecive} we find a coordinate projection $f \colon \mathbb{R}^n \to \mathbb{R}^2$ such that $f$ is injective on $\aff(S)$.
    Now
    \begin{equation*}
        f(S) = \conv(f(\chi_A),f(\chi_B),f(\chi_C),f(\chi_D)) \subseteq \mathbb{R}^2
    \end{equation*}
    is a full-dimensional square with diagonals $\conv(f(\chi_A),f(\chi_D))$ and $\conv(f(\chi_B),f(\chi_C))$ (since $f$ provides an affine isomorphism, it maps faces to faces).
    Thus $f(\chi_A) + f(\chi_D) = f(\chi_B) + f(\chi_C)$.
    From this we get $f(\chi_A) = f(\chi_B + \chi_C - \chi_D)$.
    Since $f$ is injective on $\aff(S)$ we have $\chi_A = \chi_B + \chi_C - \chi_D$ from which the result follows.
\end{proof}

\subsection{Combinatorial lemmas}
As Proposition~\ref{triangle face} suggests, we will repeatedly need to study connectivity of both order-convex sets and unions of two antichains. For later use, we observe that these notions are essentially equivalent, in the following sense.

\begin{lemma}\label{lm:extremeconn}
    Let $Q$ be a poset. Then $Q$ is connected if and only if $\max (Q)\cup \min (Q)$ is connected.
\end{lemma}
\begin{proof}
    Suppose $Q$ is connected.
    Consider any path in $Q$ between $q_1,q_2 \in \min(Q) \cup \max(Q)$.
    Without loss of generality, we may assume that this path alternates between going up in the poset and going down.
    Hence, say, 
    \begin{equation*}
        q_1 \le x_1 \ge x_2 \le x_3 \ge \cdots \le x_n \ge q_2
    \end{equation*}
    for some $x_1, \dots, x_n \in Q$
    (if $x_n \le p_2$ the following argument will work the same way).
    Let $m_1 \in \max(Q)$ be such that $x_1 \le m_1$.
    Replacing $x_1$ with $m_1$ we obtain
    \begin{equation*}
        q_1 \le m_1 \ge x_2 \le x_3 \ge \cdots \le x_n \ge q_2.
    \end{equation*}
    Let $m_2 \in \min(Q)$ be such that $x_2 \ge m_2$.
    Replacing $x_2$ with $m_2$ we obtain
    \begin{equation*}
        q_1 \le m_1 \ge m_2 \le x_3 \ge \cdots \le x_n \ge q_2.
    \end{equation*}
    We continue like this until all the $x_i$ are replaced with minimal or maximal elements.
    Hence $\min(Q) \cup \max(Q)$ is connected.
    
    Any element of $Q$ is comparable to some element in $\min(Q) \cup \max(Q)$.
    Hence, if $\min(Q) \cup \max(Q)$ is connected, so is $Q$.
\end{proof}

Connectivity of posets behaves well with respect to contractions:

\begin{lemma}\label{lm:quotients_connected}
    Let $P$ be a poset and let $J\subseteq P$ be a filter or order ideal.
    Let $S \subseteq P$ be any subset. 
    \begin{enumerate}[(a)]
        \item If $S$ is connected, then $S/J$ is connected.
        \item If both $J$ and $S/J$ are connected, then $S$ is connected.
        \item If $J$ is connected, then $S$ is connected if and only if $S/J$ is connected.
    \end{enumerate}
\end{lemma}
\begin{proof}
    \begin{enumerate}[(a)]
        \item A path $p_1 \perp p_2 \perp \cdots \perp p_k$ in $S$ induces a path  $[p_1] \perp [p_2] \perp \cdots \perp [p_k]$ in $S/J$.
        \item Let $x,y \in S$ and let $q \in J$.
        If $[x] \perp [q] \perp [y]$ then $x \perp q'$ and $q'' \perp y$ for some $q',q'' \in J$.
        Since $J$ is connected, we find a $q'$--$q''$ path inside $J$.
        Using this idea, any $[x]$--$[y]$ path in $S/J$ can be extended to an $x$--$y$ path in $S$.
        \item Follows from (a) and (b).
\end{enumerate}
\end{proof}

The following will be used implicitly in Section~\ref{subsection:triangle_ineq}.

\begin{lemma}\label{lemma:cont}
    Let $P$ be a poset, let $I$ be an order ideal and let $F$ be a filter.
    Then the following hold.
    \begin{enumerate}[(1)]
        \item If $I \neq \emptyset$ then $\min(P / I) = (\min(P) \setminus I) \cup \{ v_I \}$.
        \item If $F \neq \emptyset$ then $\max(P / F) = (\max(P) \setminus F) \cup \{ v_F \}$.
        \item If $I$ contains no maximal elements then $\max(P / I) = \max(P)$.
        \item If $F$ contains no minimal elements then $\min(P / F) = \min(P)$.
        \item If $I \cap F = \emptyset$ then $P / I / F \cong P / F / I$.
        \item Let $I = \ideal{e}$ for some $e \in P$. If $e$ is not minimal in $P$, then $P/I \cong (P - e)/(I - e)$.
        \item Let $F = \filter{e}$ for some $e \in P$. If $e$ is not maximal in $P$, then $P/F \cong (P - e)/(F - e)$.
    \end{enumerate}
\end{lemma}

\begin{proof}
  All of the above are easily verified. For an illustration, we prove part (6). 
Suppose $e$ is not minimal in $P$.
Let $f < e$.
The ground set of $P / I$ is $(P \setminus I) \cup \{ v_I \}$ and the ground set of $(P-e) / (I-e)$ is $(P-e) \setminus (I-e) \cup \{ v_{I-e} \} = (P\setminus I) \cup \{ v_{I-e} \}$.
Here $v_I$ and $v_{I-e}$ are minimal elements where for all $x \in P \setminus I$
\begin{align*}
    v_I \le x \text{ in } P / I &\iff x \ge y \le e \text{ for some } y \in P \\
    &\iff x \ge y < e \text{ or } x \ge y=e > f \\
    &\iff x \ge z \text{ for some } z \in I - e \\
    &\iff v_{I-e} \le x \text{ in } (P-e)/(I-e).
\end{align*}
Therefore $P / I \cong (P-e) / (I-e)$.
\end{proof}

As familiar from graph theory, contraction commutes with deletion, as made precise by the following lemma.

\begin{lemma}\label{lemma:del_cont}
    Let $P$ be a poset, let $J$ be either a filter or an order ideal.
    Let $S \subseteq P$ be a subset such that either $S \cap J = \emptyset$ or $J \subseteq S$.
    Then
    \begin{equation*}
        (P/J) \setminus (S/J) = (P\setminus S) / J.
    \end{equation*}
    When $J \subseteq S$ we have $(P \setminus S) \cap J = \emptyset$ and hence we may identify $(P\setminus S) / J$ with $P \setminus S$.
\end{lemma}

\begin{proof}
    Proving $(P / J) \setminus (S / J) \subseteq (P \setminus S) / J$ is straightforward.
    Suppose $[x] \in (P\setminus S) / J$ for some $x \in P \setminus S$.
    Assume towards contradiction that $[x] \in S / J$.
    Now $[x] = [x']$ for some $x' \in S$.
    We have two cases.
    
    In the first case $S \cap J = \emptyset$.
    Now $x' \not \in J$.
    Therefore $[x'] = \{ x' \}$ and hence $[x] = \{ x' \}$.
    In the second case $J \subseteq S$.
    Now $x \not \in J$ so $[x] = \{ x \}$ and hence $\{ x \} = [x']$.
    
    In both cases we have $x=x'$ which is a contradiction.
    We conclude that $[x] \in (P / J) \setminus (S / J)$ finishing the proof.
\end{proof}

The last two lemmas are not used until Section~\ref{subsection:Equality}.

\begin{lemma}\label{lemma:connected}
    Let $X \subseteq P$ be a connected subset. 
    Then $\ideal{X} \cap \filter{X}$ is connected.
\end{lemma}

\begin{proof}    
    Let $p_1$ and $p_2$ be elements from $\filter{X} \cap \ideal{X}$.
    Now $p_1 \perp x$ and $x' \perp p_2$ for some $x,x' \in X$.
    Since $X$ is connected, we find $x_1, \dots, x_n \in X$ such that
    \begin{equation*}
        p_1 \perp x \perp x_1 \perp \dots \perp x_n \perp x' \perp p_2.
    \end{equation*}
    Since $X \subseteq \filter{X} \cap \ideal{X}$, we see that $p_1$ and $p_2$ are connected by a path inside $\filter{X} \cap \ideal{X}$.
    
\end{proof}

\begin{lemma}\label{lemma:min_cup_max}
    Suppose $X \subseteq P$ is a subset of height $\le 2$.
    Let $Y = \filter{X} \cap \ideal{X}$.
    Then $\min(Y) \cup \max(Y) = X$. 
\end{lemma}

\begin{proof}
    First we claim that $\min(X) \subseteq \min(Y)$.
    Let $x \in \min(X)$.
    Now $x \in X \subseteq Y$.
    Suppose $y < x$ for some $y \in Y$.
    Now $x' \le y < x$ for some $x' \in X$, which is a contradiction since $x \in \min(X)$.
    Hence $x$ has to be minimal in $Y$, that is, $x \in \min(Y)$.
    One can prove $\max(X) \subseteq \max(Y)$ in the same way.
    Since $X$ has height $\le$ 2, we obtain
    \begin{equation*}
        X = \min(X) \cup \max(X) \subseteq \min(Y) \cup \max(Y).
    \end{equation*}
    
    Then suppose $y \in \min(Y)$.
    Now $y \in Y$ and thus  $x \le y$ for some $x \in X \subseteq Y$.
    By minimality of $y$ we cannot have $x < y$.
    Therefore $x = y$ and thus $y \in X$.
    One can prove similarly, that if $y \in \max(Y)$ then $y \in X$.
\end{proof}

\section{Triangular faces}\label{Triangles}

In this section, we will show the following theorem.

\begin{theorem}\label{thm:triangle_ineq}
Let $P$ be a poset. Then $f_2^\Delta(\Oc(P))\leq f_2^\Delta(\Cc(P))$, 
with equality holding if and only if $P$ is $\mathbf{X}$-free.
\end{theorem}

\subsection{Parametrization of triangular faces}\label{subsection:recharacterize_triangles}
We will start by giving a new combinatorial description of the triangles in both $\mathcal{O}(P)$ (Proposition~\ref{OT_Char}) and $\mathcal{C}(P)$ (Proposition~\ref{CT_Char}).
Even though these two descriptions seem more complicated than the descriptions given in Proposition~\ref{triangle face}, they will reveal a clear path towards proving Theorem~\ref{thm:triangle_ineq}.

\begin{proposition}\label{OT_Char}
    The triangles in $\mathcal{O}(P)$ are in bijection with triples $(Q,W,G)$ where 
    \begin{itemize}
        \item $Q \subseteq P$ is a connected order-convex subposet of height $\ge 2$,
        \item $W$ is an antichain of $P$ with $W \parallel Q$, and
        \item $G$ is a filter in $Q$ such that $G$ and $Q \setminus G$ are connected. 
    \end{itemize}  
\end{proposition}

We would like to remind the reader that we do not consider the empty poset to be connected.
Note that in the following proof the tuples $(Q,W,\emptyset)$ and $(Q,W,Q)$ would map to line segments rather than triangles.

\begin{proof}
    We claim that the map
    \begin{equation*}
        (Q,W,G) \longmapsto \conv(\chi_F,\chi_{\tilde G}, \chi_{H})
    \end{equation*}
    where
    \begin{align*}
        H &= \filter{Q \cup W} \\
        F &= H \setminus Q \\
        \tilde G &= F \cup G
    \end{align*}
    is a well-defined bijection from the set of triples $(Q,W,G)$ satisfying the conditions in the proposition to the set of triangles of $\mathcal{O}(P)$.

    \medbreak
    \textbf{Well-defined:}
    Given a triple $(Q,W,G)$ we need to check that $\conv(\chi_F,\chi_{\tilde G},\chi_H)$ is a triangle in $\mathcal{O}(P)$.
    We start by checking that $F,\tilde G$ and $H$ are filters of $P$.
    Clearly $H$ is a filter of $P$.
    We show that $F$ is a filter of $P$.
    Let $x \in F$ and suppose $y \ge x$.
    Now $y \in H$.
    Suppose towards contradiction that $y \in Q$.
    We have $x \ge z$ for some $z \in Q \cup W$.
    Since $Q$ is order-convex, $z \le x \le y$ and $x \not \in Q$, we cannot have $z \in Q$.
    Thus $z \in W$, but this contradicts the fact $W \parallel Q$.
    Hence $y \not\in Q$ and therefore $y \in F$.
    Thus $F$ is a filter of $P$.
    Finally, we show that $\tilde G$ is a filter of $P$.
    Let $x \in \tilde G = F \cup G$ and suppose $y \ge x$.
    If $x \in F$ then $y \in F$.
    Suppose $x \in G$.
    If $y \in Q$ then since $G \subseteq Q$ is a filter in $Q$, this implies $y \in G$.
    If $y \not \in Q$ then since $x \in Q$ we have $y \in H$ and thus $y \in H \setminus Q = F$.
    In all cases $y \in \tilde G$.
    Therefore $\tilde G$ is a filter of $P$.
    
    Next, we show that $\conv(\chi_F,\chi_{\tilde G},\chi_H)$ is a triangle in $\mathcal{O}(P)$.
    We note that $F \subseteq \tilde G \subseteq H$.
    Since $G \neq \emptyset$ we have $F \subsetneq \tilde G$.
    Since $G \subsetneq Q$ we have $F \cup G \subsetneq \filter{Q \cup W}$ and thus $\tilde G \subsetneq H$.
    We have $\tilde G \setminus F = (F \cup G) \setminus F \subseteq G$.
    If $x \in G$ then $x \in Q$ and hence $x \not\in F$.
    Hence $G \subseteq (F \cup G) \setminus F$.
    Therefore $\tilde G \setminus F = (F \cup G) \setminus F = G$ is connected. 
    Since $Q \subseteq H$ we have that 
    \begin{align*}
        H \setminus \tilde G &= H \setminus (F \cup G) \\
        &= (H \setminus F) \cap (H \setminus G) \\
        &= (H \setminus (H \setminus Q)) \cap (H  \setminus G) \\
        &= Q \cap (H \setminus G) \\
        &= Q \setminus G
    \end{align*} 
    is connected.
    Finally, $H \setminus F = H \setminus (H \setminus Q) = Q$ is connected.
    By Proposition~\ref{triangle face}(a) the set $T = \conv(\chi_F,\chi_{\tilde G},\chi_H)$ is a triangle face of $\mathcal{O}(P)$.
    Hence the map is well-defined.
    
    \medbreak
    \textbf{Injective:}
    Suppose $(Q_1,W_1,G_1)$ and $(Q_2,W_2,G_2)$ are triples satisfying the conditions of the proposition such that
    \begin{equation*}
        \conv(\chi_{F_1},\chi_{\tilde G_1},\chi_{H_1}) = \conv(\chi_{F_2},\chi_{\tilde G_2},\chi_{H_2}).
    \end{equation*}
    From this we get $\{ F_1,\tilde G_1, H_1 \} = \{ F_2, \tilde G_2, H_2 \}$.
    Since $F_1 \subsetneq \tilde G_1 \subsetneq H_1$ and $F_2 \subsetneq \tilde G_2 \subsetneq H_2$, we must have $F_1 = F_2$, $\tilde G_1 = \tilde G_2$ and $H_1 = H_2$.
    Now $F_1 = H_1 \setminus Q_1 = H_1 \setminus Q_2 = F_2$ and since both $Q_1$ and $Q_2$ are subsets of $H_1$, we get $Q_1 = Q_2$.
    We let $Q \coloneqq Q_1 = Q_2$.
    
    From $H_1 = H_2$ we obtain $\filter{Q \cup W_1} = \filter{Q \cup W_2}$.
    Let $x \in W_1$.
    Now $x \ge y \in Q \cup W_2$.
    But since $W_1 \parallel Q$, we must have $y \in W_2$.
    Similarly, $y \ge z$ for some $z \in W_1$.
    Thus $z \le y \le x$.
    Since $W_1$ is an antichain, this implies $z = x$ and hence $x = y \in W_2$.
    Therefore $W_1 \subseteq W_2$ and by symmetry $W_2 \subseteq W_1$.
    Thus $W_1 = W_2$.
    From $\tilde G_1 = \tilde G_2$ and $F_1 = F_2$ we obtain $F_1 \cup G_1 = F_1 \cup G_2$.
    Let $x \in G_1$.
    Now $x \in Q$ and hence $x \not \in F_1 = H_1 \setminus Q$.
    Therefore $x \in G_2$.
    Thus $G_1 \subseteq G_2$ and by symmetry $G_2 \subseteq G_1$.
    We have now obtained $(Q_1,W_1,G_1) = (Q_2,W_2,G_2)$.
    Hence the map is injective.
    
    \medbreak
    \textbf{Surjective:}
    Let $T = \conv(\chi_F,\chi_{\tilde G},\chi_H)$ be a triangle in $\mathcal{O}(P)$ where $F \subsetneq \tilde G \subsetneq H$.
    We let 
    \begin{align*}
        Q &= H \setminus F, \\
        W &= F \cap \min(H), \\
        G &= \tilde G \setminus F.
    \end{align*}
    Our aim is to show that $(Q,W,G)$ is in the domain of our map, \emph{i.e.}\ it satisfies the conditions of the proposition and that it maps to $\conv(\chi_F,\chi_{\tilde G},\chi_H)$.
    
    The set $H\setminus F$ is connected since $\conv(\chi_F,\chi_H)$ is an edge in $\OP$. 
    It is order-convex, because if we have $x<y<z$ with $x,z \in H \setminus F$, then $z\not\in F\Rightarrow y\not\in F$ and $x\in H\Rightarrow y\in H$, so $y\in H\setminus F$. 
    Since $F\subsetneq G\subsetneq H$, the set $H\setminus F$ has cardinality at least 2, so by connectivity it has height $\ge 2$.
    Clearly $W$ is an antichain of $P$.
    We have $F \cap \min(H) \parallel H \setminus F$ and thus $W \parallel Q$.
    We claim $G$ is a filter in $Q$.
    First we notice that $G = \tilde G \setminus F \subseteq H \setminus F = Q$.
    Suppose $x \in G$ and $y \in Q$ such that $x \le y$.
    Now $x \in \tilde G \setminus F$ and $y \in H \setminus F$.
    Since $\tilde G$ is a filter, this implies $y \in \tilde G \setminus F = G$.
    Thus $G$ is a filter in $Q$.
    Furthermore, $G \neq \emptyset$ since $F \subsetneq \tilde G$ and $G \neq Q$ since otherwise $\tilde G \setminus F = H \setminus F$ which would imply $\tilde G = H$.
    Since $\conv(\chi_F,\chi_{\tilde G},\chi_H)$ is a triangle in $\mathcal{O}(P)$ we have that $G = \tilde G \setminus F$ and $Q \setminus G = (H \setminus F) \setminus (\tilde G \setminus F) = H \setminus \tilde G$ are connected.
    Therefore the triple $(Q,W,G)$ is in the domain of our map.

    Finally, we have
    \begin{align*}
        \filter{Q \cup W} &= \filter{Q \cup (F \cap \min(H))} \\
        &= \filter{H \cap (Q \cup \min(H))}\\
        &= \filter{Q \cup \min(H)}\\
        &= \filter{Q} \cup \filter{\min(H)}\\
        &= \filter{Q} \cup H\\
        &= H
    \end{align*}
    as well as $H \setminus Q = H \setminus (H \setminus F) = F$ and $F \cup G = F \cup (\tilde G \setminus F) = \tilde G$.
    Therefore $(Q,W,G) \mapsto \conv(\chi_F,\chi_{\tilde G},\chi_H)$.
    Thus the map is surjective, finishing the proof of the proposition.
\end{proof}

The next lemma will help us give a similar characterization of the triangles in $\mathcal{C}(P)$.

\begin{lemma}\label{C_ordered}
    Let $T$ be a triangle in $\CP$. Then there are unique antichains $A, B, C$ such that $T=\conv(\chi_A,\chi_B, \chi_C)$ and $A=\min(A\cup B\cup C)$, $C=\max(A\cup B\cup C)$.
    Furthermore, in this case $A \cap C \subseteq B$.
\end{lemma}
\begin{proof}
    By Proposition~\ref{triangle face}, we have $T=\conv(\chi_A,\chi_B, \chi_C)$ for some distinct antichains $A,B,C$ such that $A\Delta B$, $A\Delta C$ and $B\Delta C$ are all connected. 
    Denote $A\cup B\cup C=S$, and after relabeling, assume without loss of generality that $A$ is the antichain among these three that has the largest intersection with $\min (S)$. 
    Suppose towards contradiction that $A\neq\min (S)$.
    Note that we cannot have $\min (S) \subsetneq A$ as every element in $S$ is comparable to some element in $\min (S)$.
    Hence we find an element $b \in \min (S) \setminus A$, and we can assume without loss of generality that $b\in B$.
    By maximality of $|A\cap \min (S)|$, we must also find some $a\in (A\setminus B) \cap \min (S)$. 
    Since a path in $A\Delta B$ must both alternate between $A$ and $B$ and alternate between up and down steps in the poset, there can not be such a path between $a$ and $b$, contradicting the connectivity of $A\Delta B$. 
    Hence, $A=\min (S)$.
    
    In the same way, one of the sets $A$, $B$, $C$ has to equal $\max (S)$.
    Suppose towards contradiction that  $A=\min (S)=\max (S)$.
    If we had some $x \in S \setminus A$ then $x$ would be neither minimal nor maximal in $S$, and hence $a < x < a'$ for some $a,a' \in A$, a contradiction.
    Hence $S = A$ and thus $B, C \subseteq A$.
    Now $A\Delta B = A \setminus B$ and $A\Delta C = A \setminus C$.
    As they are connected, they are both singletons, so $B=A\setminus\{x\}$ and $C=A\setminus\{y\}$ for some distinct $x,y\in A$. 
    Hence $B\Delta C$ is the antichain $\{x,y\}$, contradicting connectivity. 
    Therefore $\max (S)$ must equal $B$ or $C$, and after relabeling, we conclude $C=\max (S)$.
    
    Lastly, we show $A \cap C \subseteq B$.
    Let $x \in A \cap C$.
    Now $x$ is isolated in $A \cup B \cup C$.
    Suppose towards contradiction that $x \not\in B$.
    Now $x$ is isolated in $A \Delta B$ and in $C \Delta B$.
    Since both $A \Delta B$ and $C \Delta B$ are connected, we have $A \Delta B = C \Delta B = \{ x \}$.
    Now $B \setminus C = \emptyset$ and $B \setminus A = \emptyset$.
    Therefore $B \subseteq A$ and $B \subseteq C$ and hence $A \cap B = C \cap B$.
    We also have $C \setminus B = A \setminus B =  \{ x \}$.
    But now
    \begin{equation*}
        A = (A \setminus  B) \cup (A \cap B) = (C \setminus B) \cup (C \cap B) = C, 
    \end{equation*}
    a contradiction.
    We conclude $x \in B$, finishing the proof.
\end{proof}

\begin{proposition}\label{CT_Char}
    The triangles in $\mathcal{C}(P)$ are in bijection with triples $(Q,W,B)$ where
    \begin{itemize}
        \item $Q \subseteq P$ is a connected order-convex subposet of height $\ge 2$,
        \item $W \subseteq P$ is an antichain with $W \parallel Q$, and
        \item $B$ is an antichain in $Q$ such that $B \Delta \min (Q)$ and $B \Delta \max (Q)$ are connected.
    \end{itemize}
\end{proposition}

\begin{proof}
    Our aim is to show that the map
    \begin{equation*}
        (Q,W,B) \longmapsto \conv(\chi_A,\chi_{\tilde B},\chi_C),
    \end{equation*}
    where
    \begin{align*}
        A &= W \cup \min (Q),  \\
        \tilde B &= W \cup B, \\
        C &= W \cup \max (Q)
    \end{align*} 
    is a well-defined bijection from the set of triples satisfying the conditions of the proposition to the set of triangles in $\mathcal{C}(P)$.
    
    \medbreak
    \textbf{Well-defined:}
    Given a triple $(Q,W,B)$ we need to show that the set $\conv(\chi_A,\chi_{\tilde B},\chi_C)$ is a triangle in $\mathcal{C}(P)$.
    Clearly $A,\tilde B,C$ are antichains of $P$.
    Since $B \neq \min(Q)$, $B \neq \max(Q)$, $Q$ has height $\ge 2$ and $W \parallel Q$, it is easy to check that $A,\tilde B,C$ are pairwise distinct.
    Here
    \begin{equation*}
        A \Delta \tilde B = (W \cup \min (Q)) \Delta (W \cup B) = B \Delta \min (Q)
    \end{equation*}
    and 
    \begin{equation*}
        \tilde B \Delta C = (W \cup B) \Delta (W \cup \max (Q)) = B \Delta \max (Q)
    \end{equation*}
    are connected by assumption.
    Furthermore, by Lemma~\ref{lm:extremeconn}(a) 
    \begin{equation*}
        A \Delta C = (W \cup \min (Q)) \Delta (W \cup \max (Q)) = \min (Q) \Delta \max (Q) = \min (Q) \cup \max (Q)
    \end{equation*}
    is connected since $Q$ is connected.
    Therefore by Proposition~\ref{triangle face}(b), $\conv(\chi_A,\chi_{\tilde B},\chi_C)$ is a triangle in $\mathcal{C}(P)$.
    We also note that
    \begin{align*}
        \min(A \cup \tilde B \cup C) &= \min(W \cup \min (Q) \cup B \cup \max (Q)) \\
        &= W \cup \min (Q) \\
        &= A
    \end{align*}
    and 
    \begin{align*}
        \max(A \cup \tilde B \cup C) &= \max(W \cup \min (Q) \cup B \cup \max (Q)) \\
        &= W \cup \max (Q) \\
        &= C.
    \end{align*}
    
    \medbreak
    \textbf{Injective:}
    Let $(Q_1,W_1,B_1)$ and $(Q_2,W_2,B_2)$ be tuples satisfying the conditions of the proposition such that
    \begin{equation*}
        \conv(\chi_{A_1},\chi_{\tilde B_1},\chi_{C_1}) = \conv(\chi_{A_2},\chi_{\tilde B_2},\chi_{C_2}).
    \end{equation*}
    From this we obtain $\{ A_1,\tilde B_1,C_1 \} = \{ A_2,\tilde B_2,C_2 \}$.
    Therefore
    \begin{equation*}
        A_1 = \min(A_1 \cup \tilde B_1 \cup C_1) = \min(A_2 \cup \tilde B_2 \cup C_2) = A_2
    \end{equation*}
    and
    \begin{equation*}
        C_1 = \max(A_1 \cup \tilde B_1 \cup C_1) = \max(A_2 \cup \tilde B_2 \cup C_2) = C_2.
    \end{equation*}
    It follows that $\tilde B_1 = \tilde B_2$.
    We now have $W_1 \cup \min (Q_1) = W_2 \cup \min (Q_2)$ and $W_1 \cup \max (Q_1) = W_2 \cup \max (Q_2)$.
    Since $W_i \parallel Q_i$ and $\max (Q_i) \cap \min (Q_i) = \emptyset$ for $i=1,2$, this implies $\min (Q_1) = \min (Q_2)$ and $\max (Q_1) = \max (Q_2)$.
    Since $Q_1$ and $Q_2$ are order-convex, we get $Q_1 = Q_2$.
    We let $Q \coloneqq Q_1 = Q_2$.
    From $\tilde B_1 = \tilde B_2$ we get $W_1 \cup B_1 = W_2 \cup B_2$.
    Since $W_i \parallel Q$ and $B_i \subseteq Q$ for $i = 1,2$, this implies $W_1 = W_2$ and $B_1 = B_2$.
    Hence $(Q_1,W_1,B_1) = (Q_2,W_2,B_2)$.
    
    \medbreak
    \textbf{Surjective:}
    Let $\conv(\chi_A,\chi_{\tilde B},\chi_C)$ be a triangle in $\mathcal{C}(P)$ where by Lemma~\ref{C_ordered} we may assume $A = \min(A \cup \tilde B \cup C)$ and $C = \max(A \cup \tilde B \cup C)$ and $A \cap C \subseteq \tilde B$.
    We let
    \begin{align*}
        Q &= \{ x \in P \colon a \le x \le c \text{ for some } a \in A, \ c \in C \} \setminus (A \cap C), \\
        W &= A \cap C, \\
        B &= \tilde B \setminus (A \cap C).
    \end{align*}
    Our aim is to show that $(Q,W,B)$ is in the domain of our map and that it maps to $\conv(\chi_A,\chi_{\tilde B},\chi_C)$.
    
    If $Q$ contains $x$ and $z$ with $a\leq x<y<z\leq c$, then certainly $y\in Q$ as well, so $Q$ is order-convex.
    Moreover, $\max (Q) =C\setminus A$ and $\min (Q) = A\setminus C$, so $\min (Q) \cup \max (Q)  = A\Delta C$ is connected since $\conv (\chi_A,\chi_C)$ is an edge in $\mathcal{C}(P)$. 
    By Lemma~\ref{lm:extremeconn}, it follows that $Q$ is connected. 
    Since $\min (Q) \cap \max (Q) =\emptyset$, $Q$ has height $\ge 2$.
    Clearly $W$ is an antichain in $P$.
    One easily shows that $W \parallel Q$.
    Also, clearly $B$ is an antichain.
    We claim $B \subseteq Q$.
    But this is clear since $B = \tilde B \setminus (A \cap C)$ and $A = \min(A \cup \tilde B \cup C)$ and $C = \max(A \cup \tilde B \cup C)$.
    We have 
    \begin{align*}
        B \Delta \min (Q) &= (\tilde B \setminus (A\cap C)) \Delta (A \setminus C) \\
        &= \bigl( (\tilde B \setminus (A \cap C)) \cup (A \cap C) \bigr) \Delta \bigr((A \setminus C ) \cup (A \cap C)\bigr) \\
        &= \tilde B \Delta A
    \end{align*}
    which is connected since $\conv(\chi_A,\chi_{\tilde B})$ is an edge of~$\mathcal{C}(P)$.
    Similarly, we have that $B \Delta \max (Q) = \tilde B \Delta C$ is connected.
    We have now shown that $(Q,W,B)$ is in the domain of our map.
    
    Lastly, we note that 
    \begin{align*}
        W \cup \min(Q) &= (A \cap C) \cup (A \setminus C) = A, \\
        W \cup \max (Q) &= (A \cap C) \cup (C \setminus A) = C, \\
        W \cup B &= (A \cap C) \cup \bigl(\tilde B \setminus (A \cap C) \bigr) = \tilde B.
    \end{align*}
    Hence $(Q,W,B) \mapsto \conv(\chi_A,\chi_{\tilde B},\chi_C)$, finishing the proof.
\end{proof}

\subsection{Inequalities for the number of triangular faces}\label{subsection:triangle_ineq}
Using the characterizations obtained in the previous section, we will prove that the number of triangle faces of $\mathcal{O}(P)$ is at most the number of triangle faces in $\mathcal{C}(P)$.

Propositions~\ref{OT_Char} and~\ref{CT_Char} are remarkably similar in their statements. 
Using these, we can write the number of triangles in each of the polytopes as sums over connected order-convex subposets of certain combinatorial invariants of these subposets.
Let us define the invariants in question.

 \begin{definition}[Biconnected filter]\label{def:phi}
     Let $Q$ be a poset and let $X \subseteq \max (Q)$ and $Y \subseteq \min (Q)$.
     Let $G \subseteq Q$ be a filter.
     We say $G$ is \emph{biconnected} in $Q$ if $Q \setminus G$ and $G$ are connected.
     We say that $G$ is an \emph{$(X,Y)$-filter} in $Q$ if $(Q\setminus G) \cap \max (Q) = X$ and $G \cap \min (Q) = Y$.     
     
     We denote by $\varphi(Q,X,Y)$ the number of biconnected $(X, Y)$-filters in $Q$.
 \end{definition}   

\begin{definition}[Biconnected antichain]\label{def:alpha}
    Let $Q$ be a poset and let $X\subseteq \max (Q)$ and $Y\subseteq\min (Q)$.
    Let $B \subseteq Q$ be an antichain.
    We say $B$ is \emph{biconnected} in $Q$ if $B\Delta \min (Q)$ and $B\Delta \max (Q)$ are both connected. 
    We say that $B$ is an \emph{$(X,Y)$-antichain} in $Q$ if $B\cap\max (Q) = X$ and $B\cap\min (Q) = Y$. 
    
    We denote by $\alpha(Q,X,Y)$ the number of biconnected $(X,Y)$-antichains in $Q$.
\end{definition}

\begin{remark}\label{remark}
By Proposition~\ref{OT_Char} the number of triangles in $\mathcal{O}(P)$ is 
\begin{equation*}
f_2^\Delta(\Oc(P)) =   \sum_{Q} \#\{ W \colon W \text{ antichain in } P, \ W \parallel Q \} \sum_{X \subseteq \max(Q)} \sum_{Y \subseteq \min(Q)} \varphi(Q,X,Y) \label{eq:O_triangle_formula}
\end{equation*}
and by Proposition~\ref{CT_Char}, the number of triangles in $\mathcal{C}(P)$ is 
\begin{equation*}
 f_2^\Delta(\Cc(P)) =   \sum_{Q} \#\{ W \colon W \text{ antichain in } P, \ W \parallel Q \} \sum_{X \subseteq \max(Q)} \sum_{Y \subseteq \min(Q)} \alpha(Q,X,Y), \label{eq:C_triangle_formula}
\end{equation*}
where in both cases the outer sum is over connected order-convex subposets $Q \subseteq P$ of height $\ge 2$.
To finish the proof of the first part of Theorem~\ref{thm:triangle_ineq}, it remains to prove the inequality $\varphi(Q,X,Y) \le \alpha(Q,X,Y)$ for all such $Q,X,Y$. To prove the second part, we also need to show that if $P$ is not $\mathbf{X}$-free, then $P$ has a connected order-convex subposet $Q$ and $X\subseteq \max (Q)$, $Y\subseteq \min (Q)$ for which strict inequality $\varphi(Q,X,Y) < \alpha(Q,X,Y)$ holds.  
\end{remark}

\begin{remark}
  Recall that we do not consider the empty set to be connected.
  Therefore, the empty filter is never biconnected.
  However, note that the empty antichain is biconnected in $Q$ if and only if $\min(Q)$ and $\max(Q)$ are singletons.
\end{remark}

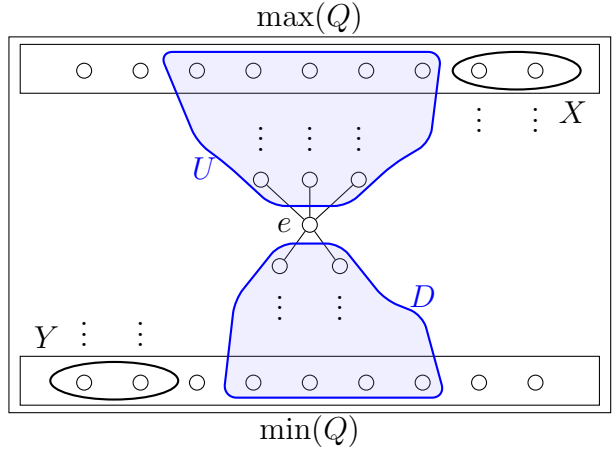
\begin{figure}
\centering
\begin{tikzpicture}[
    scale=0.5,
    every node/.style={font=\large},
    dot/.style={circle,draw,fill=white,inner sep=2.0pt},
]
\draw (0,0) rectangle (16,10);

\draw (0.3,8.5) rectangle (15.7,9.8);
\node at (8,10.8) {$\max(Q)$};

\draw (0.3,0.2) rectangle (15.7,1.5);
\node at (8,-0.8) {$\min(Q)$};

\foreach \x in {2,3.5,5,6.5,8,9.5,11,12.5,14}{
  \node[dot] at (\x,9.1) {};
}
\node at (12.5,8) {$\vdots$};
\node at (14,8) {$\vdots$};

\foreach \x in {2,3.5,5,6.5,8,9.5,11,12.5,14}{
  \node[dot] at (\x,0.8) {};
}
\node at (2,2.3) {$\vdots$};
\node at (3.5,2.3) {$\vdots$};

\draw[thick] (13.5,9.1) ellipse [x radius=1.7, y radius=0.5];
\node at (15,8) {$X$};

\draw[thick] (2.8,0.85) ellipse [x radius=1.7, y radius=0.5];
\node at (1,2) {$Y$};

\node (e) [dot,label=left:{$e$}] at (8,5) {};

\node (u1) [dot] at (6.7,6.2) {};
\node (u2) [dot] at (8.0,6.2) {};
\node (u3) [dot] at (9.3,6.2) {};
\draw (e)--(u1) (e)--(u2) (e)--(u3);

\node (d1) [dot] at (7.2,3.9) {};
\node (d2) [dot] at (8.8,3.9) {};
\draw (e)--(d1) (e)--(d2);

\filldraw[
  draw=blue, thick, rounded corners=5pt,
  fill=blue!40, fill opacity=0.15
] (4.0,9.6) -- (11.5,9.6) -- (11.2,7.2) -- (10.2,6.6) -- (9,5.5)
  -- (7,5.5) -- (5.8,6.6) -- (5,7.2) -- cycle;
\node[blue] at (4.8,6.5) {$U$};

\filldraw[
  draw=blue, thick, rounded corners=5pt,
  fill=blue!40, fill opacity=0.15
] (6.0,3) -- (7.2,4.5) -- (8.8,4.5) -- (10.0,3)
  -- (11,2.6) -- (11.6,0.4) -- (5.7,0.4) -- cycle;
\node[blue] at (11.0,3.5) {$D$};

\node at (6.7,7.5) {$\vdots$};
\node at (8.0,7.5) {$\vdots$};
\node at (9.3,7.5) {$\vdots$};
\node at (7.2,3) {$\vdots$};
\node at (8.8,3) {$\vdots$};

\end{tikzpicture}
\caption{Schematic illustration of Setup~\ref{setup}, used in Lemmas~\ref{lm:alpha-recusion}--\ref{lm:phi-recusion}.}
\label{fig:setup}
\end{figure}

We will derive recursive formulas for $\alpha(Q,X,Y)$ and $\varphi(Q,X,Y)$. 
In what follows, we will use the following setup, illustrated in  Figure~\ref{fig:setup}.

\begin{setup}\label{setup}
Let $e$ be a non-extremal element of a connected poset $Q$, let $Q'=Q\setminus \{e\}$, and let  $$D=\ideal{e}\cap Q' = \{x \in Q : x<e\}\mbox{ and }U=\filter{e}\cap Q'=\{x \in Q : x>e\}.$$ 
By non-extremality of $e$, $U$ and $D$ will then be non-empty, but not necessarily connected. We will also let $X\subseteq\max (Q)$ and $Y\subseteq\min (Q)$ be sets such that $\{e\}$, $X$ and $Y$ are pairwise parallel. 
\end{setup}

Since $e$ is parallel to $X$ and $Y$, these sets do not intersect $U$ or $D$, so we can  think of $X$ and $Y$ as sets of maximal and minimal elements in contractions of $Q$ and $Q'$ by $U$ and $D$.
By these observations, the invariants in Lemmas~\ref{lm:alpha-recusion}--\ref{lm:phi-recusion} are well-defined.
Note also that by Lemmas~\ref{lm:extremeconn}~and~\ref{lm:quotients_connected} the poset $Q'$ and all its contractions are connected and of height at least 2.

\begin{lemma}\label{lm:alpha-recusion}
  With notation as in Setup~\ref{setup},
$$\alpha(Q,X,Y)=\alpha(Q',X,Y) +\alpha(Q'/U/D,X,Y).$$ 
\end{lemma}

\begin{proof}
    Let $B$ be an antichain in $Q$ with $e \in B$.
    Now $B \cap U = \emptyset$ and $B \cap D = \emptyset$.
    Hence we may think of $B-e$ as an antichain in $Q'/U/D$ not containing $v_U$ or $v_D$.
    Thus it should be clear that we have a bijection
    \begin{align*}
        \{ \text{antichains in $Q$ containing $e$} \} &\longleftrightarrow \{ \text{antichains in $Q'/U/D$ not containing $v_U$ or $v_D$} \} \\
        B &\longleftrightarrow B-e.
    \end{align*}
    Let $B$ be an antichain in $Q$ with $e \in B$.
    Now $B \cap U = \emptyset$.
    As $\max(Q'/U) = \max(Q') \setminus U \cup \{ v_U \}$ by Lemma~\ref{lemma:cont} and $\max(Q') = \max(Q)$, we have
    \begin{align*}
        (B-e) \Delta \max(Q'/U) &= \Big( (B-e) \setminus (\max(Q) \setminus U)\Big) \cup \Big((\max(Q) \setminus U) \setminus (B-e) \Big) \cup \{ v_U \} \\
        &= \Big( (B-e) \setminus \max(Q) \Big) \cup \Big( \max(Q) \setminus (U \cup (B-e)) \Big) \cup \{ v_U \} \\
        &= \Big( (B-e) \Delta \max(Q) \Big) \setminus U \cup \{ v_U \} \\
        &= \big(B \Delta \max(Q) \big) \setminus  F \cup \{ v_U \} \\
        &\cong \big( B \Delta \max(Q) \big) / F,
    \end{align*}
    where $F \coloneqq \{e\} \cup \big(U \cap \max(Q)\big) = \{ x \in B \Delta \max(Q) \colon x \ge e \}$ is the principal filter generated by $e$ inside $B \Delta \max(Q)$.
    Dually,
    \begin{equation*}
        (B-e) \Delta \min(Q'/D) \cong \big( B \Delta \min(Q) \big) / J,
    \end{equation*}
    where $J \coloneqq  \{ e \} \cup \big( D \cap \min(Q) \big) = \{ x \in B \Delta \min(Q) \colon x \le e \}$ is the principal ideal generated by $e$ inside $B \Delta \min(Q)$.
    Since $\max(Q'/U) = \max(Q'/U/D)$ and $\min(Q'/D) = \min(Q'/U/D)$, and since $J$ and $F$ are connected, by Lemma~\ref{lm:quotients_connected} we see that $B$ is biconnected in $Q$ if and only if $B-e$ is biconnected in $Q'/U/D$.
    Furthermore, using Lemma~\ref{lemma:cont} we get
    \begin{align*}
        (B-e) \cap \max(Q'/U/D) &= (B-e) \cap \Big( \max(Q') \setminus U \cup \{ v_U \} \Big) \\
        &= B \cap \max(Q) 
    \end{align*}
    and similarly $(B-e) \cap \min(Q'/U/D) = B \cap \min(Q)$.
    With these observations we see that $B$ is a biconnected $(X,Y)$-antichain in $Q$ if and only if $B-e$ is a biconnected $(X,Y)$-antichain in $Q'/U/D$.
    By restricting the above bijection, we see that the number of biconnected $(X,Y)$-antichains in $Q$ containing $e$ is the number of biconnected $(X,Y)$-antichains in $Q'/U/D$ not containing $v_U$ or $v_D$.
    Since $U \cap X = \emptyset$ and $D \cap Y = \emptyset$ we have $v_U \not\in X$ and $v_D \not\in Y$.
    As $v_U$ and $v_D$ are maximal and minimal in $Q'/U/D$ respectively, we see that an $(X,Y)$-antichain in $Q'/U/D$ never contains $v_U$ or $v_D$.
    We conclude that the number of biconnected antichains in $Q$ containing $e$ equals $\alpha(Q'/U/D,X,Y)$.
    
    The biconnected $(X,Y)$-antichains in $Q$ not containing $e$ are exactly the biconnected $(X,Y)$-antichains in $Q'$.
    Therefore, the number of biconnected $(X,Y)$-antichains in $Q$ equals
    \begin{equation*}
        \alpha(Q,X,Y) = \alpha(Q',X,Y) + \alpha(Q'/U/D,X,Y)
    \end{equation*}
    as claimed.
\end{proof}

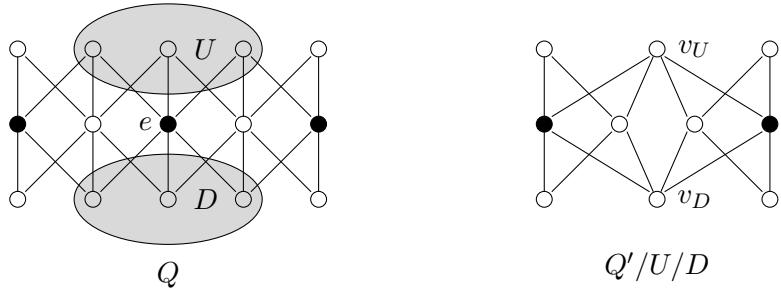
\begin{figure}[htbp]
    \centering
    \begin{tikzpicture}
\node[ellipse,
draw=black,
    fill = gray!30,
    minimum width = 2.5cm, 
    minimum height = 1.2cm] (D) at (2,0){};
\node at (2.5,0) {$D$};
\node[ellipse,
draw=black,
    fill = gray!30,
    minimum width = 2.5cm, 
    minimum height = 1.2cm] (U) at (2,2){};
    \node at (2.5,2) {$U$};
\foreach \x in {0,1,2,3,4}
{
    \draw (\x,0) circle (3pt);
    \draw (\x,2) circle (3pt);
}
\filldraw (0,1) circle (3pt);
\draw (1,1) circle (3pt);
\filldraw
(2,1) circle (3pt);

\node at (1.7,1) {$e$};

\draw (3,1) circle (3pt);
\filldraw (4,1) circle (3pt);
\foreach \x in {0,1,2,3,4}
{
    \draw (\x,0.1) -- (\x,0.9);
    \draw (\x, 1.1) -- (\x,1.9);
}
\foreach \x in {0,1,2,3}
{
    \draw (\x+0.1,0.1) -- (\x+0.9,0.9);
    \draw (\x+0.1, 1.1) -- (\x+0.9,1.9);
    \draw (\x+0.9,0.1) -- (\x+0.1,0.9);
    \draw (\x+0.9, 1.1) -- (\x+0.1,1.9);
}
\node at (2,-1) {$Q$};

\begin{scope}[xshift=7cm]
\foreach \x in {0,1.5,3}
{
    \draw (\x,0) circle (3pt);
    \draw (\x,2) circle (3pt);
}
\filldraw (0,1) circle (3pt);
\draw (1,1) circle (3pt);
\draw (2,1) circle (3pt);
\filldraw (3,1) circle (3pt);
\node at (2,2) {$v_U$};
\node at (2,0) {$v_D$};
\draw (0,0.1) -- (0,0.9);
\draw (0.1,0.1) -- (0.9,0.9);
\draw (0,1.9) -- (0,1.1);
\draw (0.1,1.9) -- (0.9,1.1);

\draw (1.4,1.9) -- (0.1,1.1);
\draw (1.45,1.9) -- (1.1,1.1);
\draw (1.55,1.9) -- (1.9,1.1);
\draw (1.6,1.9) -- (2.9,1.1);

\draw (1.4,0.1) -- (0.1,0.9);
\draw (1.45,0.1) -- (1.1,0.9);
\draw (1.55,0.1) -- (1.9,0.9);
\draw (1.6,0.1) -- (2.9,0.9);

\draw (3,0.1) -- (3,0.9);
\draw (2.9,0.1) -- (2.1,0.9);
\draw (3,1.9) -- (3,1.1);
\draw (2.9,1.9) -- (2.1,1.1);

\node at (1.5,-0.9) {$Q'/U/D$};
\end{scope}
\end{tikzpicture}
    \caption{
    On the left, a biconnected antichain (filled) in $Q$ containing the 
    element $e$. On the right, the corresponding biconnected antichain in $Q'/U/D$. Here, we have $X=Y=\emptyset$ for simplicity. 
        }\label{fig:achains'}
\end{figure}

\begin{lemma}\label{lm:alpha-supermod}
 With notation as in Setup~\ref{setup},
  we have $$\alpha(Q',X,Y)+\alpha(Q'/U/D,X,Y)\geq\alpha(Q'/U,X,Y) + \alpha(Q'/D,X,Y).$$ Moreover, if $\min(Q) \subseteq D$, $\max(Q) \subseteq U$, and $|\min(Q)|\geq 2$ and $|\max(Q)|\geq 2$, then the inequality is strict.
\end{lemma}

\begin{proof}
    Observe that all $(X,Y)$-antichains in $Q'/U$ are contained in $Q'\setminus U$, as $v_U$ is a maximal element in $Q'/U$ that is not contained in $X$. Similarly, all $(X,Y)$-antichains in $Q'/D$ are contained in $Q'\setminus D$, and all $(X,Y)$-antichains in $Q'/U/D$ are contained in $Q'\setminus (U\cup D)$. 
    
    Now, for an antichain $B$ and a poset $R$, write $t_B(R) = 1$ if $B$ is biconnected in $R$ and $0$ otherwise. 
    We can thus write
    \begin{align}
        \alpha(Q',X,Y) &= \sum_{\substack{B \subseteq Q' \\ B \text{ an $(X,Y)$-antichain in Q'}}} t_B(Q'), \label{eq:idea1}\\
        \alpha(Q'/U,X,Y) &= \sum_{\substack{B \subseteq Q' \setminus U \\ B \text{ an $(X,Y)$-antichain in  Q'}}} t_B(Q'/U), \label{eq:idea2} \\ 
        \alpha(Q'/D,X,Y) &= \sum_{\substack{B \subseteq Q' \setminus D \\ B \text{ an $(X,Y)$-antichain in  Q'}}} t_B(Q'/D), \label{eq:idea3} \\
        \alpha(Q'/U/D,X,Y) &= \sum_{\substack{B \subseteq Q' \setminus (U \cup D) \\ B \text{ an $(X,Y)$-antichain in  Q'}}} t_B(Q'/U/D). \label{eq:idea4}
    \end{align}
    Our aim is to prove $(\ref{eq:idea1}) + (\ref{eq:idea4}) \ge (\ref{eq:idea2}) + (\ref{eq:idea3})$.
    We split each sum according to whether $B$ intersects $U$ or $D$.
    In what follows, the sums are over $(X,Y)$-antichains $B \subseteq Q'$.
    Note that no antichain in $Q'$ can intersect both $U$ and $D$.
    By grouping terms together we may write
    \begin{align*}
        (\ref{eq:idea1}) + (\ref{eq:idea4}) = \sum_{B \cap D \neq \emptyset} t_B(Q') + \sum_{B \cap U \neq \emptyset} t_B(Q') + \sum_{B \cap (U \cup D) = \emptyset} \Big[ t_B(Q') + t_B(Q'/U/D) \Big] 
    \end{align*}
    and
    \begin{equation*}
        (\ref{eq:idea2}) + (\ref{eq:idea3}) = \sum_{B \cap D \neq \emptyset} t_B(Q'/U) + \sum_{B \cap U \neq \emptyset} t_B(Q'/D) + \sum_{B \cap (U \cup D) = \emptyset} \Big[ t_B(Q'/U) + t_B(Q'/D) \Big].
    \end{equation*}
    We claim that we have term-wise inequality, that is,
    \begin{enumerate}[(a)]
        \item $t_B(Q') \ge t_B(Q'/U) \text{ for all $(X,Y)$-antichains } B \subseteq Q' \setminus U, \ B \cap D \neq \emptyset$,
        \item $t_B(Q') \ge t_B(Q'/D) \text{ for all $(X,Y)$-antichains } B \subseteq Q' \setminus D, \ B \cap U \neq \emptyset$, and
        \item $t_B(Q') + t_B(Q'/U/D) \ge t_B(Q'/U) + t_B(Q'/D) \text{ for all $(X,Y)$-antichains } B$ such that $B \subseteq Q' \setminus (U \cup D)$.
    \end{enumerate}

    We start with proving (a).
    Let $B \subseteq Q' \setminus U$ such that $B \cap D \neq \emptyset$ and suppose $B$ is biconnected in $Q' / U$.
    Let $d \in B \cap D$.
    Now $B \Delta \max(Q'/U)$ and $B \Delta \min(Q' / U) = B \Delta \min(Q')$ are connected.
    Consider any path in $B \Delta \max(Q' / U)$.
    Every time this path contains a subpath $b_1 < v_U > b_2$ for some $b_1,b_2 \in B$, this means that we have $b_1 < u_1 > d < u_2 > b_2$ in $B \Delta \max(Q')$ for some $u_1,u_2 \in \max(Q') \cap U$.
    Therefore this path can be extended to a path in $B \Delta \max(Q')$.
    Therefore $B$ is biconnected in $Q'$.
    Hence we have (a). One can prove (b) in a similar way.
    
    It remains to show (c).
    Let $B$ be an $(X,Y)$-antichain in $Q'$ with $B \subseteq Q' \setminus (U \cup D)$.
    We have $B\Delta\min (Q'/U)= B\Delta\min(Q')$ and $B\Delta\max (Q'/U)$ is a contraction of $B\Delta\max(Q')$ by $U\cap\max(Q')$.
    Therefore, by Lemma~\ref{lm:quotients_connected}(a), if $B$ is biconnected in $Q'$, then it is biconnected in $Q'/U$. 
    In other words, we have $t_B(Q'/U)\geq t_B(Q')$. 
    For the same reason $t_B(Q'/U/D)\geq t_B(Q'/D)$. 
    
    If $t_B(Q'/U)= t_B(Q')$, then (c) follows. 
    It only remains to study the case where $t_B(Q'/U)> t_B(Q')$, meaning that $B$ is biconnected in $Q'/U$ but not in $Q'$. 
    In other words, this means that $B\Delta\max(Q'/U) = B \Delta \max(Q'/U/D)$ and $B \Delta\min (Q' / U) = B\Delta \min(Q') $ are both connected but $B\Delta\max(Q') = B \Delta \max(Q' / D)$ is not connected.
    By Lemma~\ref{lm:quotients_connected}, contracting the poset $B \Delta \min(Q' / U)$ with $\min(Q' / U) \cap D$ gives us the connected poset $B \Delta \min(Q'/U/D)$.
    It now follows that $B$ is biconnected in $Q'/U/D$ but not in $Q'/D$. 
    Therefore, in this case we have $t_B(Q'/U)=t_B(Q'/U/D)=1$ and $t_B(Q')=t_B(Q'/D)=0$, whence (c) follows.
    This finishes the proof of the first part of the lemma.
    
    For the second part of the lemma, we assume that $\min(Q)=\min(Q')\subseteq D$, $\max( Q)=\max (Q')\subseteq U$, and that $|\min (Q)|\geq 2$ and $|\max(Q)|\geq 2$. 
    Here $\max(Q')$ and $\min(Q')$ are not connected.
    The empty antichain $\emptyset$ is not biconnected in $Q'$ since $\emptyset \Delta \max(Q') = \max(Q')$ is not connected.
    Also, the empty antichain $\emptyset$ is not biconnected in $Q'/U$ since
    \begin{equation*}
        \emptyset \Delta \min(Q' / U) = \min(Q'/U) \stackrel{\text{Lemma \ref{lemma:cont}(4)}}{=} \min(Q')
    \end{equation*}
    is not connected, and similarly, $\emptyset$ is not biconnected in $Q' / D$ since
    \begin{equation*}
        \emptyset \Delta \max(Q' / D) = \max(Q' / D) \stackrel{\text{Lemma }\ref{lemma:cont}(3)} = \max(Q')
    \end{equation*}
    is not connected.
    However, the empty antichain $\emptyset$ is biconnected in $Q'/U/D$.
    This is because
    \begin{align*}
        \emptyset \Delta \max(Q'/U/D) &= \max(Q'/U/D) \\
        &= \max(Q' / D / U) & \mbox{[Lemma \ref{lemma:cont}(5)]} \\
        &= (\max(Q' / D) \setminus U) \cup \{ v_U \} &\mbox{[Lemma \ref{lemma:cont}(2)]} \\
        &= (\max(Q') \setminus U) \cup \{ v_U \} &\mbox{[Lemma \ref{lemma:cont}(3)]} \\
        &= \emptyset \cup \{ v_U \} &\mbox{[$\max(Q') \subseteq U$]} \\
        &= \{ v_U \}
    \end{align*}
    is connected, and similarly $\emptyset \Delta \min(Q'/U/D)$ is a singleton and hence connected.
    We therefore have
    \begin{equation*}
        t_\emptyset(Q')+t_\emptyset(Q'/U/D)=0+1>0+0=t_\emptyset(Q'/U)+t_\emptyset(Q'/D)
    \end{equation*} 
    so the inequality in (c) is strict.
    
\end{proof}

\begin{lemma}\label{lm:phi-recusion}
    With notation as in Setup~\ref{setup},
    $$\varphi(Q,X,Y)=\varphi(Q'/U,X,Y) +\varphi(Q'/D,X,Y).$$
\end{lemma}

\begin{proof}
   First, consider filters not containing $e$.
    For a filter $G \subseteq Q$ we have $e \not\in G$ if and only if $\ideal{e} \cap G = \emptyset$.
   We have a bijection
    \begin{align*}
        \Big\{ \text{filters in $Q$ not containing  $e$} \Big\} &\longleftrightarrow \Big\{ \text{filters in $Q / \ideal{e}$ not containing $v_{\ideal{e}}$ } \Big\} \\
        G &\longleftrightarrow G / \ideal{e}.
    \end{align*}
    Let $G \subseteq Q$ be a filter not containing $e$.
    Since $G \cap \ideal{e} = \emptyset$, we have $G\cong G/\ideal{e}$, so clearly $G$ is connected if and only if $G / \ideal{e}$ is connected.
    Since $\ideal{e}$ is connected, by Lemma~\ref{lm:quotients_connected}, $Q \setminus G$ is connected if and only if $(Q \setminus G) / \ideal{e}$ is connected.
    Furthermore, by Lemma~\ref{lemma:del_cont} we have $(Q \setminus G) / \ideal{e} = (Q / \ideal{e}) \setminus (G / \ideal{e})$.
    Hence we see that $G$ is biconnected in $Q$ if and only if $G / \ideal{e}$ is biconnected in $Q / \ideal{e}$. 
   
   Moreover, by Lemma~\ref{lemma:cont}, $\max(Q / \ideal{e}) = \max(Q)$ and
    \begin{align*}
        \min(Q / \ideal{e}) \cap G &= \Big( (\min(Q) \setminus \ideal{e}) \cup \{ v_{\ideal{e}} \} \Big) \cap G \\
        &= \min(Q) \cap G.
    \end{align*}
    Since $X$ and $Y$ do not intersect $\ideal{e}$, we may think of them as subsets of $\max(Q/\ideal{e})$ and $\min(Q / \ideal{e})$ respectively.
    Therefore, $G$ is an $(X,Y)$-filter in $Q$ if and only if $G / \ideal{e}$ is an $(X,Y)$-filter in $Q / \ideal{e}$.
    By restricting the above bijection we obtain that the number of biconnected $(X,Y)$-filters in $Q$ not containing $e$ equals the number of biconnected $(X,Y)$-filters in $Q/\ideal{e}$ not containing $v_{\ideal{e}}$.
   
    Since $Y \cap \ideal{e} = \emptyset$ we have $v_{\ideal{e}} \not\in Y$.
    If $G/\ideal{e}$ is an $(X,Y)$-filter in $Q/\ideal{e}$ then $v_{\ideal{e}}$ is a minimal element not in $Y = \min(Q / \ideal{e}) \cap (G / \ideal{e})$, and thus $v_{\ideal{e}} \not\in G/\ideal{e}$.
    Hence, in $Q / \ideal{e}$, the biconnected $(X,Y)$-filters not containing $v_{\ideal{e}}$ are exactly the biconnected $(X,Y)$-filters.
    Therefore the number of biconnected $(X,Y)$-filters in $Q$ not containing $e$ equals $\varphi(Q/\ideal{e},X,Y)$.
    The sets $X$ and $Y$ can be thought of as subsets of $\max(Q'/D)$ and $\min(Q'/D)$ respectively, and by Lemma~\ref{lemma:cont}, $Q / \ideal{e} \cong Q' / D$.
    Therefore the number of biconnected $(X,Y)$-filters in $Q$ not containing $e$ is
    \begin{equation*}
        \varphi(Q / \ideal{e},X,Y) = \varphi(Q' / D, X,Y).
    \end{equation*}

    Next, $(X,Y)$-filters $G\subseteq Q$ with $e\in G$ are in bijection with $(Y,X)$-filters $J\subseteq Q^{\textrm{op}}$ with $e\not \in J$ in the opposite poset, via the bijection $G\leftrightarrow J=Q\setminus G$. Indeed, for an $(X,Y)$-filter $G\subseteq Q$ we have 
    $$
    Y=\min (Q)\cap G=\min (Q)\setminus J=\max (Q^{\textrm{op}})\setminus J 
    $$ 
    and 
    $$
    X=\max (Q)\setminus G=\max (Q)\cap J=\min (Q^{\textrm{op}})\cap J.
    $$
    Moreover, by definition of biconnectivity, we see that $G$ is biconnected in $Q$ if and only if $J$ is biconnected in $Q^{\textrm{op}}$. Finally, by the first part of the proof, there are 
    \begin{align*}
        \varphi\left((Q^{\textrm{op}}-e)/\{x\in Q^\textrm{op} : x<_{Q^\textrm{op}} e\},Y,X\right) &=\varphi\left((Q^{\textrm{op}}-e)/\{x\in Q : x> e\},Y,X\right)\\
        &=\varphi\left((Q'/U)^{\textrm{op}},Y,X\right)\\&=\varphi(Q'/U,X,Y)
    \end{align*} 
    many such filters.
    
    Therefore, the number of biconnected $(X,Y)$-filters in $Q$ not containing $e$ is equal to $\varphi(Q'/D,X,Y)$, and the number of biconnected $(X,Y)$-filters in $Q$ that do contain $e$ is equal to $\varphi(Q'/U,X,Y)$.
    We conclude that the number of biconnected $(X,Y)$-filters in $Q$ equals
    \begin{equation*}
        \varphi(Q,X,Y) = \varphi(Q'/U,X,Y) + \varphi(Q'/D,X,Y)
    \end{equation*}
    as claimed.
\end{proof}

  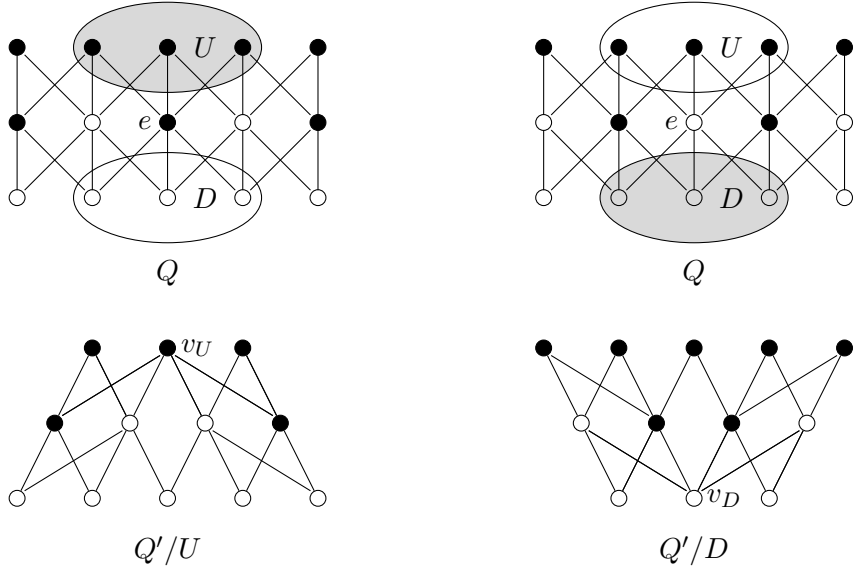
\begin{figure}[htbp]
    \centering
    \begin{tikzpicture}

\node[ellipse,
draw=black,
    minimum width = 2.5cm, 
    minimum height = 1.2cm] (D) at (2,0){};
\node at (2.5,0) {$D$};
\node[ellipse,
draw=black,
fill=gray!30,
    minimum width = 2.5cm, 
    minimum height = 1.2cm] (U) at (2,2){};
    \node at (2.5,2) {$U$};
\filldraw (0,2) circle (3pt);
\filldraw (1,2) circle (3pt);
\filldraw (2,2) circle (3pt);
\filldraw (3,2) circle (3pt);
\filldraw (4,2) circle (3pt);
\filldraw (0,1) circle (3pt);
\draw (1,1) circle (3pt);
\filldraw
(2,1) circle (3pt);
\node at (1.7,1) {$e$};
\draw (3,1) circle (3pt);
\filldraw (4,1) circle (3pt);
\draw (0,0) circle (3pt);
\draw (1,0) circle (3pt);
\draw (2,0) circle (3pt);
\draw (3,0) circle (3pt);
\draw (4,0) circle (3pt);
\node at (2,-1) {$Q$};
\foreach \x in {0,1,2,3,4}
{
    \draw (\x,0.1) -- (\x,0.9);
    \draw (\x, 1.1) -- (\x,1.9);
}
\foreach \x in {0,1,2,3}
{
    \draw (\x+0.1,0.1) -- (\x+0.9,0.9);
    \draw (\x+0.1, 1.1) -- (\x+0.9,1.9);
    \draw (\x+0.9,0.1) -- (\x+0.1,0.9);
    \draw (\x+0.9, 1.1) -- (\x+0.1,1.9);
}

\begin{scope}[yshift=-4cm]
\filldraw (1,2) circle (3pt);
\filldraw (2,2) circle (3pt);
\filldraw (3,2) circle (3pt);
\filldraw (0.5,1) circle (3pt);
\draw (1.5,1) circle (3pt);
\draw (2.5,1) circle (3pt);
\filldraw (3.5,1) circle (3pt);
 \node at (2.4,2) {$v_U$};
 \node at (2,-0.7) {$Q'/U$};
\draw (0,0) circle (3pt);
\draw (1,0) circle (3pt);
\draw (2,0) circle (3pt);
\draw (3,0) circle (3pt);
\draw (4,0) circle (3pt);
\draw (0.95,1.9) -- (0.55,1.1);
\draw (1.95,1.9) -- (1.55,1.1);
\draw (2.95,1.9) -- (2.55,1.1);
\draw (1.9,1.9) -- (0.6,1.1);
\draw (2.1,1.9) -- (3.4,1.1);
\draw (1.05,1.9) -- (1.45,1.1);
\draw (2.05,1.9) -- (2.45,1.1);
\draw (3.05,1.9) -- (3.45,1.1);
\draw (0.05,0.1) -- (0.45,0.9);
\draw (0.1,0.1) -- (1.4,0.9);
\draw (1.05,0.1) -- (1.45,0.9);
\draw (2.05,0.1) -- (2.45,0.9);
\draw (3.05,0.1) -- (3.45,0.9);
\draw (0.95,0.1) -- (0.55,0.9);
\draw (1.95,0.1) -- (1.55,0.9);
\draw (2.95,0.1) -- (2.55,0.9);
\draw (3.9,0.1) -- (2.6,0.9);
\draw (3.95,0.1) -- (3.55,0.9);
\draw (1.9,1.9) -- (0.6,1.1);
\draw (2.1,1.9) -- (3.4,1.1);
\draw (1.05,1.9) -- (1.45,1.1);
\draw (2.05,1.9) -- (2.45,1.1);
\draw (3.05,1.9) -- (3.45,1.1);
\end{scope}
\begin{scope}[xshift=7cm]
\node[ellipse,
draw=black,
    fill = gray!30,
    minimum width = 2.5cm, 
    minimum height = 1.2cm] (D) at (2,0){};
\node at (2.5,0) {$D$};
\node[ellipse,
draw=black,
    minimum width = 2.5cm, 
    minimum height = 1.2cm] (U) at (2,2){};
    \node at (2.5,2) {$U$};
\filldraw (0,2) circle (3pt);
\filldraw (1,2) circle (3pt);
\filldraw (2,2) circle (3pt);
\filldraw (3,2) circle (3pt);
\filldraw (4,2) circle (3pt);
\draw (0,1) circle (3pt);
\filldraw (1,1) circle (3pt);
\draw
(2,1) circle (3pt);
\node at (1.7,1) {$e$};
\filldraw (3,1) circle (3pt);
\draw (4,1) circle (3pt);
\draw (0,0) circle (3pt);
\draw (1,0) circle (3pt);
\draw (2,0) circle (3pt);
\draw (3,0) circle (3pt);
\draw (4,0) circle (3pt);
\node at (2,-1) {$Q$};
\foreach \x in {0,1,2,3,4}
{
    \draw (\x,0.1) -- (\x,0.9);
    \draw (\x, 1.1) -- (\x,1.9);
}
\foreach \x in {0,1,2,3}
{
    \draw (\x+0.1,0.1) -- (\x+0.9,0.9);
    \draw (\x+0.1, 1.1) -- (\x+0.9,1.9);
    \draw (\x+0.9,0.1) -- (\x+0.1,0.9);
    \draw (\x+0.9, 1.1) -- (\x+0.1,1.9);
}
\end{scope}

\begin{scope}[xshift=7cm, yshift=-4cm]
\draw (1,0) circle (3pt);
\draw (2,0) circle (3pt);
\draw (3,0) circle (3pt);
\draw (0.5,1) circle (3pt);
\filldraw (1.5,1) circle (3pt);
\filldraw (2.5,1) circle (3pt);
\draw (3.5,1) circle (3pt);
 \node at (2.4,0) {$v_D$};
 \node at (2,-0.7) {$Q'/D$};
\filldraw (0,2) circle (3pt);
\filldraw (1,2) circle (3pt);
\filldraw (2,2) circle (3pt);
\filldraw (3,2) circle (3pt);
\filldraw (4,2) circle (3pt);
\draw (0.95,0.1) -- (0.55,0.9);
\draw (1.95,0.1) -- (1.55,0.9);
\draw (2.95,0.1) -- (2.55,0.9);
\draw (1.9,0.1) -- (0.6,0.9);
\draw (2.1,0.1) -- (3.4,0.9);
\draw (1.05,0.1) -- (1.45,0.9);
\draw (2.05,0.1) -- (2.45,0.9);
\draw (3.05,0.1) -- (3.45,0.9);
\draw (0.05,1.9) -- (0.45,1.1);
\draw (0.1,1.9) -- (1.4,1.1);
\draw (1.05,1.9) -- (1.45,1.1);
\draw (2.05,1.9) -- (2.45,1.1);
\draw (3.05,1.9) -- (3.45,1.1);
\draw (0.95,1.9) -- (0.55,1.1);
\draw (1.95,1.9) -- (1.55,1.1);
\draw (2.95,1.9) -- (2.55,1.1);
\draw (3.9,1.9) -- (2.6,1.1);
\draw (3.95,1.9) -- (3.55,1.1);
\draw (1.9,0.1) -- (0.6,0.9);
\draw (2.1,0.1) -- (3.4,0.9);
\draw (1.05,0.1) -- (1.45,0.9);
\draw (2.05,0.1) -- (2.45,0.9);
\draw (3.05,0.1) -- (3.45,0.9);
\end{scope}
\end{tikzpicture}
    \caption{Biconnected (filled) filters containing $e$ (upper left) and not containing $e$ (upper right). Below are the corresponding biconnected filters in $Q'/U$ and $Q'/D$ respectively. The grey shade is used for sets that are contracted when removing the node $e$. Again, we have chosen $X=Y=\emptyset$ for simplicity.}
    \label{fig:filters}
\end{figure}

We are now ready to prove Theorem~\ref{thm:triangle_ineq}. Recall that we want to prove the inequality $f_2^\Delta(\Oc(P))\leq f_2^\Delta(\Cc(P))$,
and that equality holds if and only if $P$ is $\mathbf{X}$-free.

\begin{proof}[Proof of Theorem~\ref{thm:triangle_ineq}]
By Remark~\ref{remark}, for the first part, we only need to show that 
\begin{equation*}
  \varphi(Q,X,Y)\leq\alpha(Q,X,Y)
\end{equation*}
for all connected posets $Q$ of height at least 2, and all $X\subseteq \max (Q)$ and $Y\subseteq\min (Q)$. 

First note that, if $X \nparallel Y$, then there exist $x\in X$ and $y\in Y$ with $y<x$, since $X\cap Y\subseteq\min (Q)\cap\max (Q)  =\emptyset$ for connected posets of height at least 2. Clearly, no antichain can contain both $x$ and $y$, nor can any filter contain $y$ without also containing~$x$. 
Therefore, $\alpha(Q,X,Y) = \varphi(Q,X,Y) = 0$ if $X \nparallel Y$.

We will now assume $X\parallel Y$.
We proceed by induction over the number of non-extremal elements $e\in Q$ such that $e\parallel (X\cup Y)$. 

As a base case, assume there is no non-extremal element that is parallel to $X\cup Y$. Then the only $(X,Y)$-filter is $G=\filter{Y}\cup (\max (Q)\setminus X)$, and the only $(X,Y)$-antichain is $B=X\cup Y$. 
In such case, we have $$
    \min (G)\subseteq \min(Q)\cup\max (Q),
$$ since all non-extremal points in $Q$ are either below $X$ (and therefore not in $G$) or above $Y$ (and therefore not minimal in $G$). Noticing that 
\begin{align*}
    \min (G)\cap \min(Q)&=Y, \\
    \min (G)\cap\max (Q)&=(\max(Q)\setminus X)\setminus\filter{Y}, \\
    \max(G)&= G\cap\max(Q)=\max(Q)\setminus X,
\end{align*} we get 
\begin{align}\label{eq:filterextreme}
    \min(G) \cup \max(G)&=Y \cup (\max(Q) \setminus X)\nonumber\\
    &=\big((X \cup Y) \setminus \max(Q)\big) \cup \big(\max(Q) \setminus (X \cup Y)\big)\nonumber\\
    &=(X \cup Y) \Delta \max(Q)\nonumber\\
    &=B\Delta\max(Q).
\end{align} 

Dually, the only $(Y,X)$-filter in $Q^\textrm{op}$ is $Q \setminus G$ and the only $(Y,X)$-antichain in $Q^\textrm{op}$ is $B$.
By \eqref{eq:filterextreme},
it follows that
\begin{align}\label{eq:idealextreme}
    \min(Q \setminus G) \cup \max(Q \setminus G) &= \max_{Q^\textrm{op}}(Q \setminus G) \cup \min_{Q^\textrm{op}}(Q \setminus G) \nonumber\\
    &= B \Delta \max(Q^\textrm{op}) \nonumber\\
    &= B \Delta \min(Q).
\end{align}
where by $\max_{Q^\textrm{op}}(Q \setminus G)$ and $\min_{Q^\textrm{op}}(Q \setminus G)$ we mean the maximal and minimal elements of $Q \setminus G$ considered as a subposet of $Q^{\textrm{op}}$.

It follows from \eqref{eq:filterextreme}, \eqref{eq:idealextreme} and Lemma~\ref{lm:extremeconn} that $B$ is a biconnected antichain in $Q$ if and only if $G$ is a biconnected filter in $Q$.
Hence, $\varphi(Q,X,Y)=\alpha(Q,X,Y)\in\{0,1\}$ holds.

  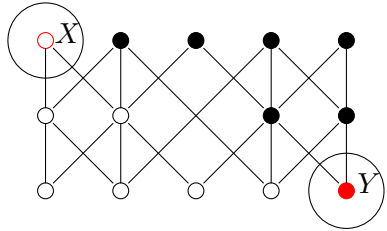
\begin{figure}[htbp]
    \centering
    \begin{tikzpicture}
\node[ellipse,
draw=black,
    minimum width = 1cm, 
    minimum height = 1cm] (Y) at (4,0){};
\node at (4.3,0.1) {$Y$};
\node[ellipse,
draw=black,
    minimum width = 1cm, 
    minimum height = 1cm] (X) at (0,2){};
\node at (0.3,2.1) {$X$};
\draw[red] (0,2) circle (3pt);
\filldraw (1,2) circle (3pt);
\filldraw (2,2) circle (3pt);
\filldraw (3,2) circle (3pt);
\filldraw (4,2) circle (3pt);
\draw (0,1) circle (3pt);
\draw (1,1) circle (3pt);
\filldraw (3,1) circle (3pt);
\filldraw (4,1) circle (3pt);
\draw (0,0) circle (3pt);
\draw (1,0) circle (3pt);
\draw (2,0) circle (3pt);
\draw (3,0) circle (3pt);
\filldraw[red] (4,0) circle (3pt);
\foreach \x in {0,1,3,4}
{
    \draw (\x,0.1) -- (\x,0.9);
    \draw (\x, 1.1) -- (\x,1.9);
}
\foreach \x in {0,2,3}
{
    \draw (\x+0.1,0.1) -- (\x+0.9,0.9);
    \draw (\x+0.9, 1.1) -- (\x+0.1,1.9);
}
\foreach \x in {0,1,3}
{
    \draw (\x+0.1, 1.1) -- (\x+0.9,1.9);
    \draw (\x+0.9,0.1) -- (\x+0.1,0.9);
}
    \draw (1.1, 1.9) -- (2.9,0.1);
    \draw (1.1,0.1) -- (2.9,1.9);
\end{tikzpicture}
    \caption{The unique $(X,Y)$-antichain (red) and the unique $(X,Y)$-filter (filled), in the base case when $Q=\filter{Y}\cup\ideal{X}\cup\max (Q)\cup\min (Q)$.}\label{fig:basecase}
\end{figure}

    For the induction step, let $e\in Q\setminus \big(\max (Q)\cup \min (Q)\cup\filter{Y}\cup\ideal{X} \big)$. Assume 
    that the inequality is proven for all triples $(R,X',Y')$ that have strictly fewer such elements than $(Q,X,Y)$ has. We define $Q'=Q\setminus\{e\}$, $D=\ideal{e}\cap Q'$ and $U=\filter{e}\cap Q'$ as in Setup~\ref{setup}.
    
    Combining this with Lemmas~\ref{lm:alpha-recusion}--\ref{lm:phi-recusion} and induction hypothesis, we get 
    \begin{align*} 
        \varphi(Q,X,Y)&=  \varphi(Q'/U,X,Y) +\varphi(Q'/D,X,Y) &\mbox{[Lemma \ref{lm:phi-recusion}]}\\&\leq \alpha(Q'/U,X,Y) +\alpha(Q'/D,X,Y)&\mbox{[Ind. hyp.]}\\ &\leq \alpha(Q',X,Y) +\alpha(Q'/U/D,X,Y)&\mbox{[Lemma \ref{lm:alpha-supermod}]}
        \\&=\alpha(Q,X,Y).&\mbox{[Lemma \ref{lm:alpha-recusion}]}
    \end{align*} 
        By the principle of induction, $\varphi(Q,X,Y)\leq\alpha(Q,X,Y)$ holds for all connected posets~$Q$ of height $\geq 2$ and all $X\subseteq \max (Q)$ and $Y\subseteq\min (Q)$.
This concludes our proof of the first part Theorem~\ref{thm:triangle_ineq}.

For the second part, we need to show that if $P$ is not $\mathbf{X}$-free, then there exists such a triple $(Q,X,Y)$ for which the inequality is strict. If $P$ is not $\mathbf{X}$-free, we find elements $d_1,d_2,e,u_1,u_2 \in P$ with $d_i<e<u_j$ for $i,j \in \{ 1,2 \}$, and $d_1 \parallel d_2$ and $u_1 \parallel u_2$. Let $Q$ be the order-convex hull of this subposet, {\em i.e.} 
\begin{equation*}
    Q=\left\{x\in P : d_i \le x \le u_j \mbox{ for some } i,j\in\{1,2\}\right\}
\end{equation*}
and let $X=Y=\emptyset$. 
Here $Q$ is a connected order-convex subposet of $P$ of height $\ge 2$.
With the selected element $e\in Q$, $U = \{ x \in Q \colon x > e \}$ and $D = \{ x \in Q \colon x < e \}$ we have $$\max (Q)=\{u_1, u_2\}\subseteq U \mbox{ and }\min (Q)=\{d_1, d_2\}\subseteq D,$$ so by Lemma~\ref{lm:alpha-supermod}, the second inequality in the chain above is strict, so $\varphi(Q,X,Y)<\alpha(Q,X,Y)$. Therefore, $f_2^\Delta(\Oc(P))<f_2^\Delta(\Cc(P))$.

Notice that if $P$ is $\mathbf{X}$-free, then $\Oc(P)$ and $\Cc(P)$ are unimodularly equivalent~\cite{HL2}, and in particular $f_2^\Delta(\Oc(P))=f_2^\Delta(\Cc(P))$. It follows that $f_2^\Delta(\Oc(P))\leq f_2^\Delta(\Cc(P))$ for every poset $P$, with equality if and only if $P$ is $\mathbf{X}$-free.
\end{proof}

\section{Square faces}\label{section:square_faces}
In this section, we will study square faces of both $\mathcal{O}(P)$ and $\mathcal{C}(P)$.
In~\cite[Theorem~7.1.1]{J}, it is shown that $f_2^\Box(\Oc(P))\leq f_2^\Box(\Cc(P))$ for any finite poset $P$. This, together with Theorem~\ref{thm:triangle_ineq}, already is enough to prove  Theorem~\ref{main} in this paper, namely that $f_2(\Oc(P))\leq f_2(\Cc(P))$ for any finite poset $P$, with equality if and only if $P$ is ${\bf X}$-free. In this section, we prove the stronger result that $f_2^\Box(\Oc(P))= f_2^\Box(\Cc(P))$ for any poset $P$.

In Section~\ref{subsection:square_of_order} we give a combinatorial description of the square faces of $\mathcal{O}(P)$, and in Section~\ref{subsection:square_chain} we give a similar description for $\mathcal{C}(P)$.
In Section~\ref{subsection:Equality} we use these descriptions to provide a bijection that shows that the number of square faces in $\mathcal{O}(P)$ and $\mathcal{C}(P)$ coincide.

\subsection{Square faces of order polytopes}\label{subsection:square_of_order}

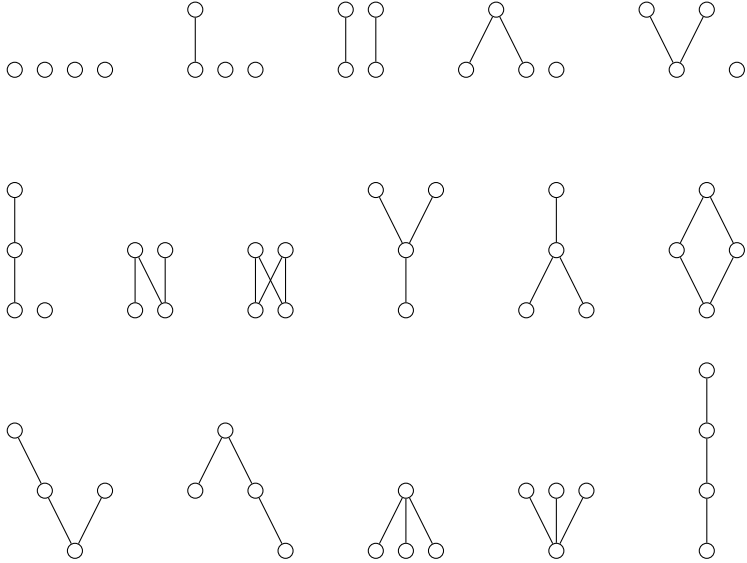
\begin{figure}[ht]
    \centering
    \begin{tikzpicture}[
        scale=0.8,
        dot/.style={draw, circle, minimum size=2mm, inner sep=0pt, fill=white},
]

        \node[dot] at (-6,0) {};
        \node[dot] at (-5.5,0) {};
        \node[dot] at (-5,0) {};
        \node[dot] at (-4.5,0) {};
        
        \node[dot] (b_l) at (-3,0) {};
        \node[dot] (b_lt) at (-3,1) {};
        \node[dot] (b_m) at (-2.5,0) {};
        \node[dot] (b_r) at (-2,0) {};
        \draw (b_l) -- (b_lt);
        
        \node[dot] (c_lt) at (-0.5,1) {};
        \node[dot] (c_lb) at (-0.5,0) {};
        \node[dot] (c_rb) at (0,0) {};
        \node[dot] (c_rt) at (0,1) {};   
        \draw (c_lt) -- (c_lb);
        \draw (c_rt) -- (c_rb);
        
        \node[dot] (d_l) at (1.5,0) {};
        \node[dot] (d_t) at (2,1) {};
        \node[dot] (d_m) at (2.5,0) {};
        \node[dot] (d_r) at (3,0) {};   
        \draw (d_l) -- (d_t);
        \draw (d_m) -- (d_t);
        
        \node[dot] (e_lt) at (4.5,1) {};
        \node[dot] (e_b) at (5,0) {};
        \node[dot] (e_mt) at (5.5,1) {};
        \node[dot] (e_r) at (6,0) {};   
        \draw (e_b) -- (e_lt);
        \draw (e_mt) -- (e_b);
        
        \node[dot] (f_lb) at (-6,-4) {};
        \node[dot] (f_lm) at (-6,-3) {};
        \node[dot] (f_lt) at (-6,-2) {};
        \node[dot] (f_r) at (-5.5,-4) {};   
        \draw (f_lb) -- (f_lm);
        \draw (f_lm) -- (f_lt);
        
        \node[dot] (g_lb) at (-4,-4) {};
        \node[dot] (g_lt) at (-4,-3) {};
        \node[dot] (g_rb) at (-3.5,-4) {};
        \node[dot] (g_rt) at (-3.5,-3) {};   
        \draw (g_lb) -- (g_lt);
        \draw (g_lt) -- (g_rb);
        \draw (g_rb) -- (g_rt);
        
        \node[dot] (h_lb) at (-1.5,-4) {};
        \node[dot] (h_lt) at (-1.5,-3) {};
        \node[dot] (h_rb) at (-2,-4) {};
        \node[dot] (h_rt) at (-2,-3) {};   

        \draw (h_lb) -- (h_lt);
        \draw (h_lt) -- (h_rb);
        \draw (h_rb) -- (h_rt);
        \draw (h_lb) -- (h_rt);
        
        \node[dot] (i_lt) at (-0,-2) {};
        \node[dot] (i_mb) at (0.5,-4) {};
        \node[dot] (i_mt) at (0.5,-3) {};
        \node[dot] (i_rt) at (1,-2) {};  
        \draw (i_lt) -- (i_mt);
        \draw (i_mt) -- (i_mb);
        \draw (i_mt) -- (i_rt);
        
        \node[dot] (j_lb) at (2.5,-4) {};
        \node[dot] (j_mb) at (3,-3) {};
        \node[dot] (j_mt) at (3,-2) {};
        \node[dot] (j_rb) at (3.5,-4) {};   
        \draw (j_lb) -- (j_mb);
        \draw (j_rb) -- (j_mb);
        \draw (j_mb) -- (j_mt);
        
        \node[dot] (k_t) at (5.5,-2) {};
        \node[dot] (k_l) at (5,-3) {};
        \node[dot] (k_r) at (6,-3) {};
        \node[dot] (k_b) at (5.5,-4) {};   
        \draw (k_b) -- (k_l);
        \draw (k_b) -- (k_r);
        \draw (k_l) -- (k_t);
        \draw (k_r) -- (k_t);
        
        \node[dot] (l_l) at (-6,-6) {};
        \node[dot] (l_lm) at (-5.5,-7) {};
        \node[dot] (l_mb) at (-5,-8) {};
        \node[dot] (l_rt) at (-4.5,-7) {};   
        \draw (l_l) -- (l_lm);
        \draw (l_lm) -- (l_mb);
        \draw (l_mb) -- (l_rt);
        
        \node[dot] (m_l) at (-3,-7) {};
        \node[dot] (m_lt) at (-2.5,-6) {};
        \node[dot] (m_rm) at (-2,-7) {};
        \node[dot] (m_rb) at (-1.5,-8) {};  
        \draw (m_l) -- (m_lt);
        \draw (m_lt) -- (m_rm);
        \draw (m_rm) -- (m_rb);
        
        \node[dot] (n_l) at (0,-8) {};
        \node[dot] (n_mb) at (0.5,-8) {};
        \node[dot] (n_mt) at (0.5,-7) {};
        \node[dot] (n_r) at (1,-8) {}; 
        \draw (n_l) -- (n_mt);
        \draw (n_mb) -- (n_mt);
        \draw (n_r) -- (n_mt);
        
        \node[dot] (o_l) at (2.5,-7) {};
        \node[dot] (o_mb) at (3,-8) {};
        \node[dot] (o_mt) at (3,-7) {};
        \node[dot] (o_r) at (3.5,-7) {};   
        \draw (o_l) -- (o_mb);
        \draw (o_mt) -- (o_mb);
        \draw (o_r) -- (o_mb);
        
        \node[dot] (p_bb) at (5.5,-8) {};
        \node[dot] (p_b) at (5.5,-7) {};
        \node[dot] (p_t) at (5.5,-6) {};
        \node[dot] (p_tt) at (5.5,-5) {};   
        \draw (p_bb) -- (p_b);
        \draw (p_b) -- (p_t);
        \draw (p_t) -- (p_tt);

    \end{tikzpicture}
    \caption{All posets with four elements}
    \label{4elements}
\end{figure}

We start by giving an explicit description of the square faces of $\mathcal{O}(P)$.
We approach this by seeing how four filters in $P$ corresponding to vertices of a square face of $\mathcal{O}(P)$ can be contained in each other as sets.

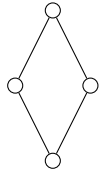
\begin{figure}
    \begin{center}
        \begin{tikzpicture}[
        scale=1,
        dot/.style={draw, circle, minimum size=2mm, inner sep=0pt, fill=white},
]
            \node[dot] (top) at (0,2) {};
            \node[dot] (a) at (-0.5,1) {};
            \node[dot] (b) at (0.5,1) {};
            \node[dot] (bottom) at (0,0) {};
            \draw (bottom) -- (a);
            \draw (bottom) -- (b);
            \draw (a) -- (top);
            \draw (b) -- (top);
        \end{tikzpicture}
    \end{center}
    \caption{The only way the filters corresponding to a square face in $\mathcal{O}(P)$ can be contained in each other.}
    \label{fig:O(P)_square_poset}
\end{figure}

\begin{lemma}\label{lemma:poset_of_filters}
Let $P$ be a poset, and let $A,B,C,D$ be pairwise distinct filters of $P$.
If $\conv(\chi_A, \chi_B, \chi_C, \chi_D)$ is a square face of $\mathcal{O}(P)$, then the poset $(\{A,B,C,D\}, \subseteq)$ is isomorphic to the poset shown in Figure~\ref{fig:O(P)_square_poset}. 
\end{lemma}

\begin{proof}
    When we order the filters $A,B,C,D$ by inclusion, then by Proposition \ref{edge} every filter is comparable to at least two other filters.
    Since $\dim(\chi_A,\chi_B,\chi_C,\chi_D) = 2$, by Lemma \ref{Lemma:princip3} there cannot be principal order ideals or principal filters of size 3 in the poset $(\{ A,B,C,D \},\subseteq)$.
    Thus by going through all 4-element posets (shown in Figure~\ref{4elements}), the only possible 4-element poset is the one shown in Figure~\ref{fig:O(P)_square_poset}.
\end{proof}

\begin{theorem}\label{SQorder}
    The square faces of $\mathcal{O}(P)$ are exactly the convex hulls of the form
    \begin{equation*}
        \conv(\chi_{F_1 \cap F_2}, \chi_{F_1}, \chi_{F_2}, \chi_{F_1 \cup F_2}),
    \end{equation*}
    where $F_1$ and $F_2$ are filters such that $F_1 \setminus F_2$ and $F_2 \setminus F_1$ are connected.
\end{theorem}

\begin{proof}
    Suppose $A,B,C,D \subseteq P$ are filters such that $\conv(\chi_A,\chi_B,\chi_C,\chi_D)$ is a square face of $\mathcal{O}(P)$.
    By Lemma~\ref{lemma:poset_of_filters}, the poset $(\{ A,B,C,D \},\subseteq)$ is isomorphic to the poset shown in Figure~\ref{fig:O(P)_square_poset}.
    By renaming the filters if necessary, we may assume the bottom element is $A$ and the top element is $D$.
    Here $A \subseteq B \cap C$.
    If we had a proper inclusion $A \subsetneq B \cap C$ then we would find $i \in (B \cap C) \setminus A$ and thus $i \in (B \cap C \cap D) \setminus A$.
    But this is not possible by Corollary~\ref{cor:onlyone}.
    Hence $A = B \cap C$.
    One can show similarly that we must have $D = B \cup C$.
    Using Proposition \ref{edge}, we see that the edges of the square must be $\conv(\chi_A,\chi_B)$, $\conv(\chi_A,\chi_C)$, $\conv(\chi_B,\chi_D)$ and $\conv(\chi_C,\chi_D)$.
    Thus $D \setminus B = (C \cup B) \setminus B = C \setminus B$ and $D \setminus C = (B \cup C) \setminus C = B \setminus C$ are connected.
    Hence by writing $F_1 = B$ and $F_2 = C$ we see that $\conv(\chi_A,\chi_B,\chi_C,\chi_D)$ is as claimed.
    
    We now prove that if $F_1$ and $F_2$ are as in the statement of the theorem, then
    \begin{equation*}
        \mathcal{F} \coloneqq \conv(\chi_{F_1 \cap F_2}, \chi_{F_1}, \chi_{F_2}, \chi_{F_1 \cup F})
    \end{equation*}
    is a square face of $\mathcal{O}(P)$.
    We prove this with induction on $|P \setminus (F_1 \cup F_2)| + |F_1 \cap F_2|$.
    
    In the base case we have $F_1 \cup F_2 = P$ and $F_1 \cap F_2 = \emptyset$.
    Now $P = F_1 \sqcup F_2$.
    Note that $F_1 \setminus F_2 = F_1$ and $F_2 \setminus F_1 = F_2$ are connected.
    Therefore by Proposition~\ref{edge}, $\conv(\chi_{\emptyset},\chi_{F_1})$ is an edge of $\mathcal{O}(F_1)$ and $\conv(\chi_\emptyset,\chi_{F_2})$ is an edge of $\mathcal{O}(F_2)$.
    Since $F_1 \cap F_2 = \emptyset$ and $F_1$ and $F_2$ are filters, we must have $F_1 \parallel F_2$.
    Thus $\mathcal{O}(F_1 \sqcup F_2) = \mathcal{O}(F_1) \times \mathcal{O}(F_2)$.
    Hence
    \begin{align*}
        \mathcal{F} = \conv(\chi_\emptyset,\chi_{F_1},\chi_{F_2},\chi_{F_1 \sqcup F_2}) = \conv(\chi_\emptyset,\chi_{F_1}) \times \conv(\chi_\emptyset,\chi_{F_2})
    \end{align*}
    is a 2-dimensional face of $\mathcal{O}(F_1) \times \mathcal{O}(F_2) = \mathcal{O}(F_1 \sqcup F_2) = \mathcal{O}(P)$ with 4 vertices, \emph{i.e.}\ a square face.
    
    Suppose $P \neq F_1 \cup F_2$.
    Now there exists a minimal element $p \in P \setminus (F_1 \cup F_2)$.
    Hence $\mathcal{F}$ is contained in the facet $\{ x_p = 0 \}$.
    We may consider $F_1$ and $F_2$ as filters of $P-p$.
    By induction assumption
    \begin{equation*}
        \conv(\chi_{F_1 \cap F_2},\chi_{F_1},\chi_{F_2},\chi_{F_1 \cup F_2}) \subseteq \mathbb{R}^{P-p}
    \end{equation*}
    is a square face of $\mathcal{O}(P-p) \cong \mathcal{O}(P) \cap \{ x_p = 0 \}$.
    Thus $\mathcal{F}$ is a square face of $\mathcal{O}(P)$.
    
    If $F_1 \cap F_2 \neq \emptyset$ then there exists  a maximal element $p \in F_1 \cap F_2$.
    Now $\mathcal{F}$ is contained in the facet $\{ x_p = 1 \}$.
    Here $F_1 - p$ and $F_2 - p$ are filters of $P-p$ such that $(F_1 - p) \setminus (F_2 - p) = F_1 \setminus F_2$ and $(F_2 - p) \setminus (F_1 - p) = F_2 \setminus F_1$ are connected.
    By induction assumption
    \begin{equation*}
        \conv(\chi_{(F_1-p) \cap (F_2-p)} , \chi_{F_1 - p}, \chi_{F_2 - p}, \chi_{(F_1 -p) \cup (F_2 - p)} )\subseteq \mathbb{R}^{P -p}
    \end{equation*}
    is a square face of $\mathcal{O}(P-p) \cong \mathcal{O}(P) \cap \{ x_p = 1 \}$.
    Therefore $\mathcal{F}$ is a square face of $\mathcal{O}(P)$.
\end{proof}

\begin{corollary}\label{cor:O(P)_bijection}
    The square faces of $\mathcal{O}(P)$ are in bijection with the set
    \begin{equation*}
        S_{\mathcal{O}(P)} \coloneqq
        \left\{ \{ F_1,F_2 \} \;\middle|\; 
        \begin{array}{l}
            F_1 \text{ and } F_2 \text{ filters such that} \\
            F_1 \setminus F_2 \text{ and } F_2 \setminus F_1 \text{ are connected}
        \end{array}
        \right\}.
    \end{equation*}
\end{corollary}

\begin{proof}
    Theorem \ref{SQorder} shows that the map from the set $S_{\mathcal{O}(P)}$ to the set of square faces of $\mathcal{O}(P)$ defined by 
    \begin{equation*}
        \{ F_1,F_2 \} \mapsto \conv(\chi_{F_1 \cap F_2}, \chi_{F_1}, \chi_{F_2},\chi_{F_1 \cup F_2})
    \end{equation*}
    is well-defined and surjective.
    It remains to show injectivity.
    Let $\{ F_1,F_2 \}$ and $\{ F_1',F_2' \}$ be elements from the set  $S_{\mathcal{O}(P)}$.
    Let 
    \begin{align*}
        Q_1 &= \conv(\chi_{F_1 \cap F_2}, \chi_{F_1},\chi_{F_2},\chi_{F_1 \cup F_2}), \\
        Q_2 &= \conv(\chi_{F'_1 \cap F'_2}, \chi_{F'_1},\chi_{F'_2},\chi_{F'_1 \cup F'_2})
    \end{align*}
    and suppose $Q_1 = Q_2$.
    We claim that $F_1 \cap F_2 = F_1' \cap F_2'$.
    Let $p \in F_1 \cap F_2$.
    Now every vertex of $Q_1$ is on the hyperplane $\{ x_p = 1 \}$.
    Hence so is every vertex of $Q_2$.
    In particular, $p \in F_1' \cap F_2'$.
    Thus $F_1 \cap F_2 \subseteq F_1' \cap F_2'$.
    Similarly, $F_1' \cap F_2' \subseteq F_1 \cap F_2$.
    We notice that we can write
    \begin{align*}
        Q_1 &= \conv(\chi_{\emptyset},\chi_{F_1 \setminus F_2}, \chi_{F_2 \setminus F_1}, \chi_{F_1 \Delta F_2}) + \chi_{F_1 \cap F_2}, \\
        Q_2 &= \conv(\chi_{\emptyset},\chi_{F'_1 \setminus F'_2}, \chi_{F'_2 \setminus F'_1}, \chi_{F'_1 \Delta F'_2}) + \chi_{F'_1 \cap F'_2}.
    \end{align*}
    Since $F_1 \cap F_2 = F_1' \cap F_2'$ we get
    \begin{equation*}
        \conv(\chi_{\emptyset},\chi_{F_1 \setminus F_2}, \chi_{F_2 \setminus F_1}, \chi_{F_1 \Delta F_2}) = \conv(\chi_{\emptyset},\chi_{F'_1 \setminus F'_2}, \chi_{F'_2 \setminus F'_1}, \chi_{F'_1 \Delta F'_2}) 
    \end{equation*} 
    and thus
    \begin{equation*}
        \{ F_1 \setminus F_2, F_2 \setminus F_1, F_1 \Delta F_2\} = \{ F_1' \setminus F_2', F_2' \setminus F_1', F_1' \Delta F_2' \}.
    \end{equation*}
    It follows that $F_1 \Delta F_2 = F_1' \Delta F_2'$.
    We have 
    \begin{equation*}
        F_1 = (F_1 \cap F_2) \cup (F_1 \setminus F_2) \subseteq (F_1' \cap F_2') \cup (F_1' \setminus F_2') \cup (F_2' \setminus F_1').
    \end{equation*}
    Since $F_1'$ and $F_2'$ are filters, we must have $(F_1' \setminus F_2') \parallel (F_2' \setminus F_1')$.
    Since $F_1 \setminus F_2$ is connected, it has to either be contained in $F_1' \setminus F_2'$ or in $F_2' \setminus F_1'$.
    In either case, we have that $F_1$ is contained in either $F_1'$ or $F_2'$.
    Similarly, $F_2$ is contained in either $F_1'$ or $F_2'$.
    By symmetry, $F_1'$ is contained in either $F_1$ or $F_2$, and likewise, $F_2'$ is contained in either $F_1$ or $F_2$.
    Since we cannot have $F_i \subseteq F_j$ or $F'_i \subseteq F'_j$ for $i \neq j$ (since all the differences are non-empty), this implies $\{ F_ 1,F_2 \} = \{ F_1',F_2' \}$.
\end{proof}

\subsection{Square faces of chain polytopes}\label{subsection:square_chain}

Next, we give the characterization of square faces of the chain polytope $\Cc(P)$ in terms of antichains of $P$.
Similarly to the case of $\mathcal{O}(P)$, we will study the square faces of $\mathcal{C}(P)$ by seeing how four antichains corresponding to the vertices of a square face can be contained in each other.

\begin{lemma}\label{antichain}
Let $P$ be a poset, and let $A$ and $B$ be distinct antichains of $P$.
If $A \Delta B$ is connected in $P$, then, up to renaming, 
\begin{itemize}
    \item [(i)] $A \subsetneq B$ and $|B \setminus A| = 1$, or  
    \item[(ii)] $A \nsubseteq B$, and $a < b$ holds whenever $a \in A$ and $b \in B$ are comparable and distinct.
\end{itemize}
\end{lemma}

\begin{proof}
In the proof of \cite[Theorem 2.1]{HLSS}, the following claim is obtained (up to renaming): If $A \Delta B$ is connected, then $A$ and $B$ satisfy either (i) $A \subsetneq B$, or (ii) $A \nsubseteq B$ and $a < b$ holds whenever $a \in A$ and $b \in B$ are comparable.
In case (i), we have $A \Delta B = B \setminus A$.
Since $B \setminus A$ is an antichain and connected, it follows that $|B \setminus A|=1$.  
\end{proof}

Let $P$ be a poset, and let $\Ac(P)$ be the set of antichains of $P$. 
We let $(\Ac(P), \subseteq)$ denote the poset of $\Ac(P)$ ordered by inclusion.

\begin{proposition}\label{cover of A(P)}
Let $A_i$, $A_j \in \Ac(P)$ be distinct antichains.
Then $A_i$ is covered by $A_j$ in the poset $(\Ac(P), \subseteq)$ if and only if $A_j = A_i \sqcup \{p\}$ for some $p \in P$. 
In this case $\conv(\chi_{A_i}, \chi_{A_j})$ is an edge of $\Cc(P)$.
\end{proposition}

\begin{proof}
The first claim is clear.
Suppose $A_j = A_i \sqcup \{p\}$ for some $p \in P$. Then $A_i \Delta A_j = \{ p \}$, which is connected.
Thus, $\conv(\chi_{A_i},\chi_{A_j})$ is an edge of $\Cc(P)$ by Proposition \ref{edge}.
\end{proof}

\begin{proposition}\label{cover relation}
Let $P$ be a poset, and let $A_1$, $A_2$, $A_3$, and $A_4$ be pairwise distinct antichains of $P$.
Suppose that $\Fc = \conv(\chi_{A_1}, \chi_{A_2}, \chi_{A_3}, \chi_{A_4})$ is a square face of $\Cc(P)$. 
If $A_i$ is covered by $A_j$ in the poset $(\{A_1, A_2, A_3, A_4\}, \subseteq)$ for some $i, j \in [4]$, then $A_i$ is also covered by $A_j$ in the poset $(\Ac(P), \subseteq)$.
\end{proposition}

\begin{proof}
    Suppose towards contradiction that $A_i$ is not covered by $A_j$ in $\mathcal{A}(P)$.
    Now $A_i \subsetneq A_j$ and $|A_j \setminus A_i| \ge 2$.
    Therefore $A_i \Delta A_j = A_j \setminus A_i$ is an antichain of size $\ge 2$ and thus not connected.
    Hence $\conv(\chi_{A_i},\chi_{A_j})$ is not an edge of $\mathcal{C}(P)$.
    Let $\{ k,\ell \} = [4] \setminus \{ i,j \}$.
    Here the diagonals of the square face $\mathcal{F}$ are $\conv(\chi_{A_i},\chi_{A_j})$ and $\conv(\chi_{A_k},\chi_{A_\ell})$.
    By Corollary~\ref{cor:diagonals} we have $A_i \cap A_j = A_k \cap A_\ell$  and $A_i \cup A_j = A_k \cup A_\ell$.
    Since $A_i \subseteq A_j$ we obtain
    \begin{equation*}
        A_i = A_i \cap A_j = A_k \cap A_\ell \subseteq A_k \subseteq A_k \cup A_\ell = A_i \cup A_j = A_j,
    \end{equation*}
    that is, $A_i \subseteq A_k \subseteq A_j$.
    This contradicts the assumption that $A_j$ covers $A_i$ in the poset $(\{ A_1,A_2,A_3,A_4 \},\subseteq)$.
\end{proof}

From Proposition \ref{cover of A(P)} and Proposition \ref{cover relation}, we obtain the following corollary.
This implies that if $\{\chi_{A_1}, \chi_{A_2}, \chi_{A_3}, \chi_{A_4}\}$ forms a square face of $\Cc(P)$, then the cover relations in the poset $(\{A_1, A_2, A_3, A_4\}, \subseteq)$ correspond to edges of $\Cc(P)$ that are of type~(i) in Lemma~\ref{antichain}. 

\begin{corollary}\label{cover cor}
Let $P$ be a poset, and let $A_1$, $A_2$, $A_3$, and $A_4$ be pairwise distinct antichains of $P$.
Suppose that $\conv(\chi_{A_1}, \chi_{A_2}, \chi_{A_3}, \chi_{A_4})$ is a square face of $\Cc(P)$. 
If $A_i$ is covered by $A_j$ in the poset $(\{A_1, A_2, A_3, A_4\}, \subseteq)$ for some $i, j \in [4]$, then there exists some $p \in P$ such that $A_j = A_i \sqcup \{p\}$.
In particular, $\conv(\chi_{A_i}, \chi_{A_j})$ is an edge of $\Cc(P)$.
\end{corollary}

Using the above arguments, we determine the Hasse diagram of $(\{A, B, C, D\}, \subseteq)$ that corresponds to a square face of $\Cc(P)$.

\begin{proposition}\label{ABCD}
Let $P$ be a poset, and let $A$, $B$, $C$, and $D$ be pairwise distinct antichains of $P$.
If $\Fc = \conv(\chi_{A}, \chi_{B}, \chi_{C}, \chi_{D})$ is a square face of $\Cc(P)$, then the
poset
\begin{equation*}
    (\{A, B, C, D\}, \subseteq)
\end{equation*}
is isomorphic to one of the posets shown in Figure~\ref{fig:abc}.

\medskip
\begin{figure}[htbp]
\centering
\begin{tikzpicture}[
    scale=1,
    dot/.style={draw, circle, minimum size=2mm, inner sep=0pt, fill=white},
]
    \node[dot] (a_top) at (-3,2) {};
    \node[dot] (a_left) at (-3.5,1) {};
    \node[dot] (a_right) at (-2.5,1) {};
    \node[dot] (a_bottom) at (-3,0) {};
    \draw (a_bottom) -- (a_left);
    \draw (a_bottom) -- (a_right);
    \draw (a_left) -- (a_top);
    \draw (a_right) -- (a_top);
    \node at (-3,-0.7) {{\rm(a)}};

    \node[dot] (b_ltop) at (0.5,1) {};
    \node[dot] (b_lbottom) at (0.5,0) {};
    \node[dot] (b_rbottom) at (1.5,0) {};
     \node[dot] (b_rtop) at (1.5,1) {};
    \draw (b_ltop) -- (b_lbottom);
    \draw (b_rtop) -- (b_rbottom);
    \node at (1,-0.7) {{\rm(b)}};

    \node[dot] (c1) at (4,0) {};
    \node[dot] (c2) at (4.5,0) {};
    \node[dot] (c3) at (5,0) {};
    \node[dot] (c4) at (5.5,0) {};
    \node at (4.75,-0.7) {\rm{(c)}};s

\end{tikzpicture}
\caption{The only possible ways the antichains corresponding to a square face of $\mathcal{C}(P)$ can be contained in each other.}
\label{fig:abc}
\end{figure}
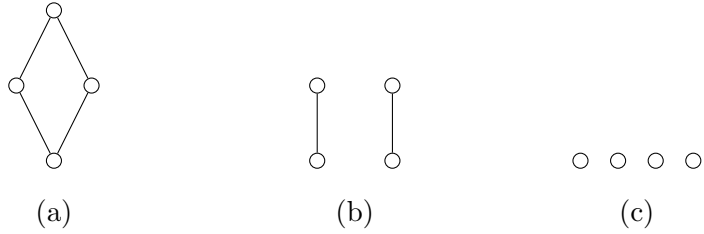
\end{proposition}

\begin{figure}
    \centering
    \begin{tikzpicture}[
    scale=1,
        dot/.style={draw, circle, minimum size=2mm, inner sep=0pt, fill=white},
]
        \node[dot] (d_top) at (-3.5,1) {};
        \node[dot] (d_bottom) at (-3.5,0) {};
        \node[dot] (d_middle) at (-3,0) {};
        \node[dot] (d_right) at (-2.5,0) {};
        \draw (d_bottom) -- (d_top);
        \node at (-3,-0.7) {{\rm(d)}};
        
        \node[dot] (e_left) at (0.5,0) {};
        \node[dot] (e_middle) at (1,0) {};
        \node[dot] (e_right) at (1.5,0) {};
        \node[dot] (e_top) at (1,1) {};
        \draw (e_top) -- (e_left);
        \draw (e_top) -- (e_middle);
        \draw (e_top) -- (e_right);
        \node at (1,-0.7) {{\rm(e)}};
        
        \node[dot] (f_left) at (4.25,1) {};
        \node[dot] (f_middle) at (4.75,1) {};
        \node[dot] (f_right) at (5.25,1) {};
        \node[dot] (f_bottom) at (4.75,0) {};
        \draw (f_bottom) -- (f_left);
        \draw (f_bottom) -- (f_middle);
        \draw (f_bottom) -- (f_right);
        \node at (4.75,-0.7) {\rm{(f)}};s
    \end{tikzpicture}
    \caption{Posets considered in the proof of Proposition \ref{ABCD}.}
    \label{fig:def}
\end{figure}
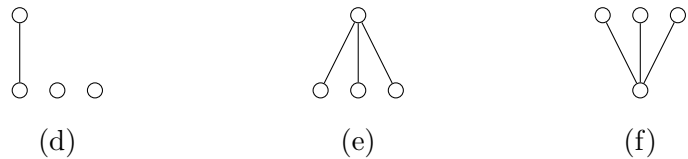

\begin{proof}
Since $\dim\Fc = 2$, by Lemma \ref{Lemma:princip3} the poset $(\{A, B, C, D\}, \subseteq)$ cannot have principal ideals or principal filters of size 3.
By going through all 4-element posets (shown in Figure~\ref{4elements}), the poset $(\{ A,B,C,D \},\subseteq)$ has to be isomorphic to either one of the posets shown in Figure~\ref{fig:abc} or to one of the posets shown in Figure~\ref{fig:def}.
Our aim is to show that none of the cases shown in Figure~\ref{fig:def} are possible.

When the poset $(\{A, B, C, D\}, \subseteq)$ is isomorphic to either Figure~\ref{fig:def} (e) or (f), three edges are incident to a single vertex in $\mathcal{F}$ by Corollary~\ref{cover cor}.  
That is, $\dim \Fc = 3$.
This contradicts the assumption that $\mathcal{F}$ is a square face.

Suppose that the poset $(\{A, B, C, D\}, \subseteq)$ is isomorphic to Figure~\ref{fig:def} (d), where $\emptyset \neq A \subsetneq B$.
By Corollary~\ref{cover cor}, we know that $B = A \sqcup \{p\}$ for some $p\in P$.
Let the diagonals of $\Fc$ be the segments $\conv(\chi_A,\chi_D)$ and $\conv(\chi_B,\chi_C)$.
Then, $\conv(\chi_A, \chi_C)$ is an edge of $\Cc(P)$ and $\conv(\chi_B, \chi_C)$ is not an edge of $\Cc(P)$.
That is, $A \Delta C$ is connected in $P$, while $B \Delta C$ is disconnected. 
We claim that $p \not\in C$.
If not then $p$ belongs to both antichains $B$ and $C$, and hence $p$ is an isolated point in $B \cup C = A\cup B \cup C$.
In particular, since $p \in C \setminus A$, $p$ is isolated in $A \Delta C$.
Since $A \not\subseteq C$ and $C \not\subseteq A$ we have $|A \Delta C| \ge 2$.
This contradicts the fact that $A \Delta C$ is connected. 
Thus, we have $p \notin C$.

Next we claim that $p$ is not comparable to anything in $C \setminus B$.
We have 
\begin{equation*}
    B \Delta C = (A \Delta C) \sqcup \{ p \},
\end{equation*}
where $B \Delta C$ is not connected and $A \Delta C$ is connected.
If $p$ was comparable to something in $C \setminus B$, then $p$ would be comparable to something in $B \Delta C$ and hence comparable to something in $A \Delta C$.
But since $A \Delta C$ is connected this would mean that also $B \Delta C$ is connected, which is not the case.
Thus $p$ is not comparable to anything in $C \setminus B$.

The element $p$ is also not comparable to anything in $C \cap B$ since otherwise $p$ is comparable to something in $B$ but this is not possible since $p \in B$ and $B$ is an antichain.
Hence $p$ is not comparable to anything in $C$ and thus $C \sqcup \{ p \}$ is an antichain.

Now, let $h \colon \mathbb{R}^P \to \mathbb{R}$ be a linear map and let $\alpha \in \mathbb{R}$ be such that the equation $h(x) \le \alpha$ defines $\Fc$ as a face of $\mathcal{C}(P)$.
We have $h(\chi_A) = h(\chi_B) = \alpha$ and thus $h(\chi_{\{ p \}}) = 0$.
Hence $h(\chi_{C \sqcup \{ p \}}) = \alpha$.
This implies that $C \sqcup \{ p \} \in \{ A,B,C,D \}$.
Thus we need to have $C \sqcup \{ p \} = B$ or $C \sqcup \{ p \} = D$ and consequently $C \subseteq B$ or $C \subseteq D$.
This contradicts the fact that the poset $(\{ A,B,C,D \},\subseteq )$ is of the form (d).
\end{proof}

\begin{theorem}\label{thm:square_faces_of_C(P)}
    Let $P$ be a poset.
    The square faces of $\mathcal{C}(P)$ are exactly the sets of the form
    \begin{equation*}
        \conv(\chi_{Q_1 \cup R_1 \cup S}, \chi_{Q_1 \cup R_2 \cup S}, \chi_{Q_2 \cup R_1 \cup S}, \chi_{Q_2 \cup R_2 \cup S}),
    \end{equation*}
    where $S$ is an antichain of $P$, and $Q$ and $R$ are connected subposets of $P$ satisfying 
    \begin{itemize}
        \item $Q = Q_1 \sqcup Q_2$ where $Q_1$ and $Q_2$ are antichains of $P$,
        \item $R = R_1 \sqcup R_2$ where $R_1$ and $R_2$ are antichains of $P$, and
        \item $Q \parallel R$, $Q \parallel S$ and  $R \parallel S$.
    \end{itemize}
\end{theorem}

\begin{figure}[htbp]
\centering

\begin{tikzpicture}[
  scale=0.5,
  dot/.style={circle,draw,fill=white,inner sep=2pt}
]
\fill[orange!15]
(0,0) rectangle (3.4,6.5);   
\fill[orange!15]
(4,0) rectangle (7.4,6.5);   
\fill[orange!15]
(8,0) rectangle (11.4,6.5);  

\node (q1) [dot] at (0.9,1.2) {};
\node (q2) [dot] at (1.7,1.2) {};
\node (q3) [dot] at (2.5,1.2) {};
\node (q4) [dot] at (0.9,5.3) {};
\node (q5) [dot] at (1.7,5.3) {};
\node (q6) [dot] at (2.5,5.3) {};

\draw (q1) -- (q4);
\draw (q2) -- (q5);
\draw (q3) -- (q6);
\draw (q1) -- (q5);
\draw (q2) -- (q6);

\node (r1) [dot] at (8.9,1.2) {};
\node (r2) [dot] at (9.7,1.2) {};
\node (r3) [dot] at (10.5,1.2) {};
\node (r4) [dot] at (8.9,5.3) {};
\node (r5) [dot] at (9.7,5.3) {};
\node (r6) [dot] at (10.5,5.3) {};

\draw (r1) -- (r4);
\draw (r2) -- (r5);
\draw (r3) -- (r6);
\draw (r1) -- (r5);
\draw (r2) -- (r6);

\node (s1) [dot] at (4.9,3.25) {};
\node (s2) [dot] at (5.8,3.25) {};
\node (s3) [dot] at (6.7,3.25) {};

\coordinate (C) at (5.7,3.25);

\draw[thick] (5.7,2) 
ellipse [x radius=6.5, y radius=1.6];

\draw[thick,green!80!black] (5.7,4.5) 
ellipse [x radius=6.5, y radius=1.6];

\draw[thick,red,rotate around={22:(C)}]
(C) ellipse [x radius=6.5, y radius=1.6];

\draw[thick,blue,rotate around={-22:(C)}]
(C) ellipse [x radius=6.5, y radius=1.6];

\node[font=\Large] at (1.7,-0.7) {$Q$};
\node[font=\Large] at (5.7,-0.7) {$S$};
\node[font=\Large] at (9.7,-0.7) {$R$};

\end{tikzpicture}
\caption{Antichains corresponding to the vertices of a square face of $\mathcal{C}(P)$.}
\label{fig:QRS}
\end{figure}
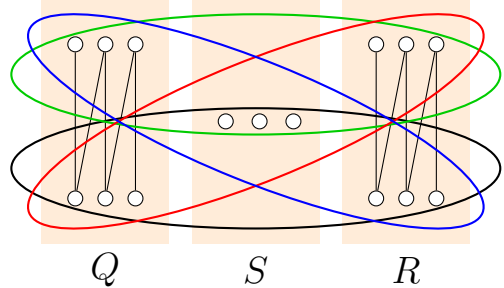

A connected subposet that is a disjoint union of two antichains is either of height 2 or a singleton.
Thus, either $Q = \min(Q) \sqcup \max(Q)$ or $Q = \{ p \} \sqcup \emptyset$.
Hence the sets $Q_1$ and $Q_2$ are uniquely determined (up to renaming) by the set $Q$, and similarly $R_1$ and $R_2$ are determined by $R$.
To illustrate Theorem~\ref{thm:square_faces_of_C(P)}, Figure~\ref{fig:QRS} shows the case when $Q$ and $R$ have height 2.

\begin{proof}
    Let $A$, $B$, $C$, and $D$ be pairwise distinct antichains of $P$.
    Suppose that $\mathcal{F} = \conv(\chi_A,\chi_B,\chi_C,\chi_D)$ is a square face of $\Cc(P)$, and let the diagonals of $\mathcal{F}$ be the segments $\conv(\chi_A,\chi_D)$ and $\conv(\chi_B,\chi_C)$.
    We aim to show that
    \begin{equation*}
        \{ A,B,C,D \} = \{Q_1 \cup R_1 \cup S,\; Q_1 \cup R_2 \cup S,\; Q_2 \cup R_1 \cup S,\; Q_2 \cup R_2 \cup S\}
    \end{equation*}
    for some suitable $Q, R$ and $S$.
    By Proposition~\ref{ABCD} the poset $(\{ A,B,C,D \}, \subseteq)$ is isomorphic to one of the posets shown in Figure~\ref{fig:abc}.
    We split into cases accordingly.
    In all cases, we choose $S = A \cap B \cap C \cap D$.
    
    \medskip
    \textbf{Case 1:} Suppose $(\{ A,B,C,D \}, \subseteq)$ is isomorphic to the poset in Figure~\ref{fig:abc}(a).
    Now by Corollary~\ref{cover cor} we may assume that the bottom element in $(\{ A,B,C,D \},\subseteq)$ is $A$  and the top element is $D$.
    Thus $S = A = A \cap B \cap C \cap D$.
    By Corollary~\ref{cover cor} we have $B = A \sqcup \{ p \}$, $C = A \sqcup \{ q \}$ and $D = A \sqcup \{ p\} \sqcup \{ q \}$ for some $p,q \in P$.
    Since $D$ is an antichain, we have $p \parallel q$, $\{p\} \parallel A$ and $\{q\} \parallel A$.
    Hence we may choose $Q = \{ q \}$, $R = \{ p \}$ with $Q_1 = \emptyset$, $Q_2 = \{ q \}$ and $R_1 = \emptyset$, $R_2 = \{ p \}$.
    
    \medskip
    \textbf{Case 2:} Suppose $(\{ A,B,C,D \},\subseteq)$ is isomorphic to the poset in Figure~\ref{fig:abc}(b).
    We may assume that $\emptyset \neq A \subseteq B$ and $\emptyset \neq  C \subseteq D$.
    By Corollary~\ref{cover cor} we have $B = A \sqcup \{ p \}$ and $D = C \sqcup \{ q \}$ for some $p,q \in P$.
    Here $p \in B \setminus A$.
    By Corollary~\ref{cor:onlyone}, $p$ cannot belong the exactly 1 or 3 of the sets $A,B,C,D$.
    If $p \in C$ then $p \in (B \cap C \cap D) \setminus A$ which is not possible.
    Thus $p \in D \setminus C$ and hence $p = q$.
    Since $\conv(\chi_A,\chi_C)$ is an edge, $A \Delta C$ is connected by Proposition~\ref{edge}.
    As $A \not\subseteq C$ and $C \not\subseteq A$ we have $|A \Delta C| \ge 2$.
    Because $A \Delta C$ is connected and cannot have height $> 2$ it must have height 2.
    We let $Q = A \Delta C$ and $R = \{ p \}$.
    Here $Q = Q_1 \sqcup Q_2$ where $Q_1 \coloneqq \min(A \Delta C)$ and $Q_2 \coloneqq \max(A \Delta C)$.
    Also $R = R_1 \sqcup R_2$ where $R_1 \coloneqq \emptyset$ and $R_2 \coloneqq \{ p \}$.
    Clearly $Q \parallel S$ and $R \parallel S$.
    Since $A \subseteq B$ and $p \in B$, $p$ cannot be comparable to any element in $A$.
    Similarly, because $C \subseteq D$ and $p \in D$, $p$ cannot be comparable to any element in $C$.
    Thus $Q \parallel R$.
    Since $A$ and $C$ satisfy the condition in Lemma \ref{antichain} (ii), either $\max(A \Delta C) \subseteq A$ and $\min(A \Delta C) \subseteq C$, or $\max(A \Delta C) \subseteq C$ and $\min(A \Delta C) \subseteq A$. 
    Let $(i,j) \in \{ (1,2),(2,1) \}$ be such that $Q_i \subseteq A$ and $Q_j \subseteq C$.
    We claim that     
    \begin{align}
        A &= Q_i \cup R_1 \cup S,   \label{eq:squareface_case2_1'} \\
        B &= Q_i \cup R_2 \cup S,   \label{eq:squareface_case2_2'} \\
        C &= Q_j \cup R_1 \cup S,  \label{eq:squareface_case2_3'} \\
        D &= Q_j \cup R_2 \cup S. \label{eq:squareface_case2_4'}
    \end{align}
    The inclusion $\supseteq$ in all cases is clear.
    Consider (\ref{eq:squareface_case2_1'}).
    Let $a \in A$.
    If $a \in A \cap C$ then $a \in A \cap B \cap C \cap D = S$.
    So suppose $a \not\in A \cap C$.
    Now $a \in A \Delta C$.
    Since $A \Delta C$ has height 2 and $Q_i \subseteq A$ we need to have $a \in Q_i$.
    Hence we have the inclusion $\subseteq$ in (\ref{eq:squareface_case2_1'}).
    Since $B = A \sqcup \{ p \}$ the inclusion $\subseteq$ in (\ref{eq:squareface_case2_2'}) is clear.
    The inclusion $\subseteq$ in (\ref{eq:squareface_case2_3'}) and (\ref{eq:squareface_case2_4'}) are proven in the same way.
    
    \medskip
    \textbf{Case 3:}
    Suppose $(\{A,B,C,D\} ,\subseteq)$ is isomorphic to the poset in Figure~\ref{fig:abc}(c).
    By Corollary~\ref{cor:onlyone}, no element of $P$ can belong to exactly 1 or 3 of the sets $A,B,C,D$.
    Thus
    $$
    A \Delta C = B \Delta D \quad \mbox{and} \quad 
    A \Delta B = C \Delta D.
    $$
    Note that both of these sets have cardinality $\ge 2$, are connected, cannot have height $> 2$, and thus have height 2.
    By Corollary~\ref{cor:diagonals} we have $A \cap D = B \cap C$ and $A \cup D = B \cup C$.

    We claim that
    \begin{equation}
        (A \Delta C) \cap (C \Delta D) = \emptyset. \label{eq:empty_delta}
    \end{equation}
    Indeed, if an element $x$ is in the set in (\ref{eq:empty_delta}) then either $x \in (A \cap D) \setminus C$ or $x \in C \setminus (A \cup D)$.
    But the facts $A \cap D = B \cap C$ and $A \cup D = B \cup C$ give rise to a contradiction in both cases.
    Therefore (\ref{eq:empty_delta}) holds.

    We also claim that 
    \begin{equation}
        (A \Delta C) \parallel (C \Delta D) \label{eq:delta_parallel}.
    \end{equation}
    Suppose towards contradiction that $p \perp q$ for some $p \in A \Delta C$ and $q \in C \Delta D$.
    Note that $p \neq q$ by (\ref{eq:empty_delta}).
    We split into four cases depending on whether $p \in A \setminus C$ or $p \in C \setminus A$ and whether $q \in C \setminus D$ or $q \in D \setminus C$. 
    As an example, we consider the case $p \in A \setminus C$ and $q \in C \setminus D$.
    In this case, $p \in A \Delta C = B \Delta D$ so either $p \in B$ or $p \in D$.
    Since $A \cap D = B \cap C$ and $p \in A \setminus C$, we must have $p \in B$.
    Similarly, $q \in C \Delta D = A \Delta B$ and thus $q$ is either in $A$ or $B$.
    But this is not possible since $p \perp q$, $p \neq q$, $p \in A \cap B$, and $A$ and $B$ are antichains.
    Therefore the case $p \in A \setminus C$ and $q \in C \setminus D$ is not possible.
    Other cases can be proven to be impossible in a similar way.
    We therefore conclude that (\ref{eq:delta_parallel}) holds.

    Let $Q = A \Delta C$ and $R = C \Delta D$.
    We saw above that $Q \parallel R$.
    Again, clearly $Q \parallel S$ and $R \parallel S$.
    Here, $Q$ and $R$ are connected and are of height 2.
    We let $Q_1 = \min(A \Delta C)$, $Q_2 = \max(A \Delta C)$, $R_1 = \min(C \Delta D)$, $R_2 = \max(C \Delta D)$.
    Here $Q_1,Q_2,R_1,R_2$ are antichains such that $Q = Q_1 \sqcup Q_2$ and $R = R_1 \sqcup R_2$.
    By Lemma \ref{antichain}(ii), we have
    either $\max(A \Delta C) \subseteq A$ and $\min(A \Delta C) \subseteq C$, or $\max(A \Delta C) \subseteq C$ and $\min(A \Delta C) \subseteq A$.
    Similarly for $C \Delta D$.
    Let $(i,j), (k,\ell) \in \{ (1,2),(2,1) \}$ 
    be such that 
    \begin{equation*}
        Q_i \subseteq A, \quad Q_j \subseteq C, \quad R_k \subseteq C, \quad R_\ell \subseteq D.
    \end{equation*}
    We claim that 
    \begin{align}
        A &= Q_i \cup R_k \cup S,  \label{eq:aqrs_1} \\
        B &= Q_i \cup R_\ell \cup S,\label{eq:aqrs_2} \\
        C &= Q_j \cup R_k \cup S,  \label{eq:cqrs_1} \\
        D &= Q_j \cup R_\ell \cup S \label{eq:cqrs_2}.
    \end{align}
    We prove (\ref{eq:aqrs_1}).
    
    For the inclusion $\subseteq$, let $x \in A$.
    As $A \cup D = B \cup C$, we have $x \in B \setminus C$, $x \in C \setminus B$ or $x \in B \cap C$.
    First, if $x \in B \cap C$ then $x \in A \cap B \cap C$, and since no element of $P$ can belong to exactly 3 of the sets $A,B,C,D$, this implies $x \in A \cap B \cap C \cap D = S$.
    Then suppose $x \in B \setminus C$.
    Now $x \in A \setminus C \subseteq A \Delta C = Q = Q_i \sqcup Q_j$, and as $Q_j \subseteq C$ we have $x \not\in Q_j$ so $x \in Q_i$.
    Lastly, suppose $x \in C \setminus B$.
    Then $x \in A \cap C$ and as $x \in D$ would imply $x \in B$, we must have $x \not\in D$.
    Hence $x \in C \Delta D = R = R_k \sqcup R_\ell$.
    Since $R_\ell \subseteq D$ we have $x \in R_k$.
    We have now shown the inclusion $\subseteq $ in (\ref{eq:aqrs_1}).
    For the inclusion $\supseteq$ it suffices to prove $R_k \subseteq A$.
    Let $x \in R_k \subseteq C$.
    Now $x \in R = C \Delta D = A \Delta B$.
    If $x \in A$ we are done, and if $x \in B$ then $x \in B \cap C = A \cap D$ so $x \in A$.
    We therefore have (\ref{eq:aqrs_1}).
    One can prove (\ref{eq:aqrs_2}), (\ref{eq:cqrs_1}) and (\ref{eq:cqrs_2}) in a similar way.

    We have now seen that each square face is as claimed.

    \medskip
    Suppose $Q,R,S$ satisfy the condition of the theorem.
    We aim to prove that
    \begin{equation*}
        \mathcal{F} = \conv(\chi_{Q_1 \cup R_1 \cup S}, \chi_{Q_1 \cup R_2 \cup S}, \chi_{Q_2 \cup R_1 \cup S}, \chi_{Q_2 \cup R_2 \cup S})
    \end{equation*}
    is a square face of $\mathcal{C}(P)$.
    We prove this with induction on $|P \setminus (Q \cup R \cup S)| + |S|$.
    
    In the base case we have $P = Q \cup R \cup S$ and $S = \emptyset$.
    Now $P = Q \sqcup R$.
    Note that  $Q_1 \Delta Q_2 = Q$.
    Since $Q$ is connected and $Q_1,Q_2$ are antichains of $Q$, by Proposition \ref{edge} the line segment $\conv(\chi_{Q_1},\chi_{Q_2})$ is an edge of $\mathcal{C}(Q)$.
    Similarly, $\conv(\chi_{R_1},\chi_{R_2})$ is an edge of $\mathcal{C}(R)$.
    Since $Q \parallel R$, we have $\mathcal{C}(Q \sqcup R) = \mathcal{C}(Q) \times \mathcal{C}(R)$.
    Thus
    \begin{align*}
       \mathcal{F} &= \conv(\chi_{Q_1 \sqcup R_1}, \chi_{Q_1 \sqcup R_2}, \chi_{Q_2 \sqcup R_1}, \chi_{Q_2 \sqcup R_2}) \\
        &= \conv(\chi_{Q_1},\chi_{Q_2}) \times \conv(\chi_{R_1},\chi_{R_2})
    \end{align*}
    is a $2$-dimensional face of $\mathcal{C}(Q) \times \mathcal{C}(R) = \mathcal{C}(Q \sqcup R) = \mathcal{C}(P)$ with $4$ vertices, \emph{i.e.}\ a square face.
    
    If there exists $p \in P \setminus (Q \cup R \cup S)$, then $\mathcal{F}$ is on the facet $\{ x_p = 0 \}$.
    We may consider $S$ as an antichain of $P-p$ and $Q$ and $R$ as connected subposets of $P-p$.
    By induction assumption
    \begin{equation*}
        \conv(\chi_{Q_1 \cup R_1 \cup S}, \chi_{Q_1 \cup R_2 \cup S}, \chi_{Q_2 \cup R_1 \cup S}, \chi_{Q_2 \cup R_2 \cup S}) \subseteq \mathbb{R}^{P - p}
    \end{equation*}
    is a square face of $\mathcal{C}(P-p) \cong \mathcal{C}(P) \cap \{ x_p = 0 \}$ and hence $\mathcal{F}$ is a square face of $\mathcal{C}(P)$.
    
    Then suppose $P = Q \cup R \cup S$.
    Suppose there exists $p \in S$. 
    Now $\{ p \}$ is a maximal chain since $Q \parallel S$ and $R \parallel S$ and $S$ is an antichain.
    Hence $\mathcal{F}$ is on the facet $\{ x_p = 1 \}$.
    Here $p \not \in Q \cup R$ since $Q \parallel S$ and $R \parallel S$.
    We may thus consider $Q$ and $R$ as connected subposets of $P-p$ and $S-p$ as an antichain of $P-p$.
    We also have $Q \parallel (S - p)$ and $R \parallel (S-p)$ in the poset $P-p$.
    By induction assumption we have that
    \begin{equation*}
        \conv(\chi_{Q_1 \sqcup R_1 \sqcup (S - p)},\chi_{Q_1 \sqcup R_2 \sqcup (S - p)},\chi_{Q_2 \sqcup R_1 \sqcup (S - p)},\chi_{Q_2 \sqcup R_2 \sqcup (S - p)}) \subseteq \mathbb{R}^{P-p}
    \end{equation*}
    is a square face of $\mathcal{C}(P-p) \cong \mathcal{C}(P) \cap \{ x_p = 1 \}$ and hence $\mathcal{F}$ is a square face of $\mathcal{C}(P)$.
\end{proof}

\begin{corollary}\label{cor:C(P)_bijection}
    The square faces of $\mathcal{C}(P)$ are in bijection with the following set, which consists of ordered pairs of two subposets and an antichain of $P$:
    \begin{equation*}
        S_{\mathcal{C}(P)} \coloneqq 
        \left\{ (\{Q,R\},S)  \;\middle|\;
        \begin{array}{l}
            Q = Q_1 \sqcup Q_2 \text{ connected, $Q_i$ antichains }  \\
            R = R_1 \sqcup R_2 \text{ connected, $R_i$ antichains } \\
            S \text{ antichain } \\
            Q \parallel S, \ Q \parallel R, \ R \parallel S
        \end{array}
        \right\}.
    \end{equation*}
\end{corollary}

Given sets $Q,R$ and $S$ satisfying the conditions of Theorem~\ref{thm:square_faces_of_C(P)}, we can construct the square face $\conv(\chi_{Q_1 \cup R_1 \cup S}, \chi_{Q_1 \cup R_2 \cup S}, \chi_{Q_2 \cup R_1 \cup S}, \chi_{Q_2 \cup R_2 \cup S})$ of $\mathcal{C}(P)$.
Recall that the sets $Q_1$ and $Q_2$ are uniquely determined (up to renaming) by the set $Q$.
Similarly for $R_1$ and $R_2$.
Here $Q$ and $R$ play symmetric roles and hence in order to have a bijection we do not want to differentiate between them.
However, the set $S$ cannot be swapped with either of the sets $Q$ or $R$ without changing the associated convex hull.
Furthermore, the set $S$ cannot be distinguished given the triple $\{ Q,R,S \}$.
For example, it is possible that $Q$, $R$ and $S$ are all singletons.
Thus, in order to recognize which set plays the role of $S$, we need the set $S_{\mathcal{C}(P)}$ to consist of ordered pairs $(\{ Q,R \},S)$ instead of unordered triples $\{ Q,R,S \}$.

\begin{proof}
    Theorem~\ref{thm:square_faces_of_C(P)} shows that the map from the set $S_{\mathcal{C}(P)}$ to the square faces of $\mathcal{C}(P)$ defined by
    \begin{equation*}
        (\{ Q,R \},S) \mapsto \conv(\chi_{Q_1 \cup R_1 \cup S}, \chi_{Q_1 \cup R_2 \cup S}, \chi_{Q_2 \cup R_1 \cup S}, \chi_{Q_2 \cup R_2 \cup S})
    \end{equation*}
    is well-defined and surjective.
    It remains to show injectivity.
    Suppose $(\{ Q,R \},S)$ and $(\{ Q',R' \},S')$ are elements from $S_{\mathcal{C}(P)}$.
    Let
    \begin{align*}
        \mathcal{F}_1 \coloneqq&\conv(\chi_{Q_1 \cup R_1 \cup S}, \chi_{Q_1 \cup R_2 \cup S}, \chi_{Q_2 \cup R_1 \cup S}, \chi_{Q_2 \cup R_2 \cup S}),  \\
        \mathcal{F}_2 \coloneqq& \conv(\chi_{Q'_1 \cup R'_1 \cup S'}, \chi_{Q'_1 \cup R'_2 \cup S'}, \chi_{Q'_2 \cup R'_1 \cup S'}, \chi_{Q'_2 \cup R'_2 \cup S'})
    \end{align*}
    and suppose $\mathcal{F}_1 = \mathcal{F}_2$.
    Our aim is to show $\{ Q,R \} = \{ Q',R' \}$ and $S = S'$.
    Suppose towards contradiction that $S \neq S'$.
    Say $p \in S \setminus S'$.
    Now every vertex of $\mathcal{F}_1$ is on the hyperplane $\{ x_p = 1 \}$ and hence so is every vertex of $\mathcal{F}_2$.
    This implies that all the sets $Q_i' \cup R_j'$, $i,j = 1,2$, contain the element $p$.
    But since $Q'_1 \cap Q_2' = R_1' \cap R_2' = Q' \cap R' = \emptyset$, we have
    \begin{equation*}
        p \in (Q_1' \cup R_1') \cap (Q_1' \cup R_2') \cap (Q_2' \cup R_1') \cap (Q_2' \cup R_2') = \emptyset,
    \end{equation*}
    a contradiction.
    Therefore $S = S'$.
    Note that we can write
    \begin{align*}
        \mathcal{F}_1 &= \conv(\chi_{Q_1 \cup R_1}, \chi_{Q_1 \cup R_2}, \chi_{Q_2 \cup R_1}, \chi_{Q_2 \cup R_2}) + \chi_S, \\
        \mathcal{F}_2 &= \conv(\chi_{Q'_1 \cup R'_1}, \chi_{Q'_1 \cup R'_2}, \chi_{Q'_2 \cup R'_1}, \chi_{Q'_2 \cup R'_2}) + \chi_{S'}.
    \end{align*}
    Since $S = S'$ we get
    \begin{equation*}
        \conv(\chi_{Q_1 \cup R_1}, \chi_{Q_1 \cup R_2}, \chi_{Q_2 \cup R_1}, \chi_{Q_2 \cup R_2}) = \conv(\chi_{Q'_1 \cup R'_1}, \chi_{Q'_1 \cup R'_2}, \chi_{Q'_2 \cup R'_1}, \chi_{Q'_2 \cup R'_2})
    \end{equation*}
    and therefore
    \begin{equation*}
        \{ Q_1 \cup R_1, Q_1 \cup R_2, Q_2 \cup R_1, Q_2 \cup R_2 \} = \{ Q'_1 \cup R'_1, Q'_1 \cup R'_2, Q'_2 \cup R'_1, Q'_2 \cup R'_2 \}.
    \end{equation*}
    For $Q_1$ we find some $i,j \in \{ 1,2 \}$ such that $Q_1 \subseteq Q_i' \cup R_j'$ and for $Q_2$ we find some $k,\ell \in \{ 1,2 \}$ such  that $Q_2 \subseteq Q_k' \cup R_\ell'$.
    Therefore
    \begin{equation*}
        Q \subseteq Q_i' \cup R_j' \cup Q_k' \cup R_\ell'.
    \end{equation*}
    Since $Q$ is connected and $Q' \parallel R'$, the set $Q$ has to be contained in either $Q_i' \cup Q_k'$ or in $R_j' \cup R_\ell'$, and hence in either $Q'$ or $R'$.
    Similarly, $R$ is contained in either $Q'$ or $R'$.
    By symmetry, $Q'$ is contained in either $Q$ or $R$ and similarly $R'$ is contained in either $Q$ or $R$.
    Since $Q$ and $R$ are not contained in each other, and since $Q'$ and $R'$ are not contained in each other, we must have $\{ Q,R \} = \{ Q',R' \}$.
\end{proof}

\subsection{Equality of the number of square faces}\label{subsection:Equality}
In this section, we prove that the number of square faces in $\mathcal{O}(P)$ and $\mathcal{C}(P)$ always coincide.
In other words, we prove that $f_2^\Box(\Oc(P))=f_2^\Box(\Cc(P))$. 
This is a strengthening of \cite[Theorem 7.1.1]{J}, where it is shown that $f_2^\Box(\Oc(P))\leq f_2^\Box(\Cc(P))$.
The maps used in the proof of Theorem~\ref{thm:square_eq} were inspired by the map constructed in  \cite[Section 7.1]{J}. There, the description of square faces of $\Cc(P)$ was incomplete, wherefore the author was unable to show surjectivity of the map analogous to our map $\Phi$ below.

We will obtain the equality $f_2^\Box(\mathcal{O}(P)) = f_2^\Box(\mathcal{C}(P))$ by constructing an explicit bijection between the sets from Corollary~\ref{cor:O(P)_bijection} and Corollary~\ref{cor:C(P)_bijection} that parametrize the square faces of each polytope.

\begin{theorem}\label{thm:square_eq}
    For any poset $P$, $f_2^\Box(\Oc(P))=f_2^\Box(\Cc(P))$.
\end{theorem}

\begin{figure}
    \centering
    \begin{tikzpicture}[scale=0.6]
        
        \node[draw, ellipse, minimum width=0.9cm, minimum height=0.7cm] (Q1) at (0,0) {$Q_1$};
        \node[draw, ellipse, minimum width=0.9cm, minimum height=0.7cm] (Q2) at (0,4) {$Q_2$};
        \node[draw, ellipse, minimum width=0.9cm, minimum height=0.7cm] (Q3) at (3,4) {};
        \node[draw, ellipse, minimum width=0.9cm, minimum height=0.7cm](S) at (5,2) {$S$};
        \node[draw, ellipse, minimum width=0.9cm, minimum height=0.7cm] (Q5) at (7,4) {};
        \node[draw, ellipse, minimum width=0.9cm, minimum height=0.7cm] (R1) at (10,0) {$R_1$};
        \node[draw, ellipse, minimum width=0.9cm, minimum height=0.7cm] (R2) at (10,4) {$R_2$};
        \node[draw, ellipse, minimum width=0.9cm, minimum height=0.7cm] (Q8) at (5,7) {};
        
        \draw[-] (Q1) -- (Q2);
        \draw[-] (Q1) -- (Q3);
        \draw[-] (Q2) -- (Q8);
        \draw[-] (Q3) -- (Q8);
        \draw[-] (S) -- (Q8);
        \draw[-] (Q5) -- (Q8);
        \draw[-] (R1) -- (R2);
        \draw[-] (R2) -- (Q8);
        \draw[-] (R1) -- (Q5);
        
        \draw[thick, blue] plot [smooth cycle] coordinates {(-0.5,-1) (6,1) (8,4) (5,8) (-1.5,4.5)};
        
        \draw[thick, red, opacity=0.5] plot [smooth cycle] coordinates {(10.5,-1) (4,1) (2,4) (5,8) (11.5,4.5)};
        
        \node[color=blue] at (0,7) {$F_1$};
        \node[color=red] at (10,7) {$F_2$};
        
    \end{tikzpicture}
    \caption{The sets involved in the maps $\Phi$ and $\Psi$.}
    \label{fig:schematic}
\end{figure}
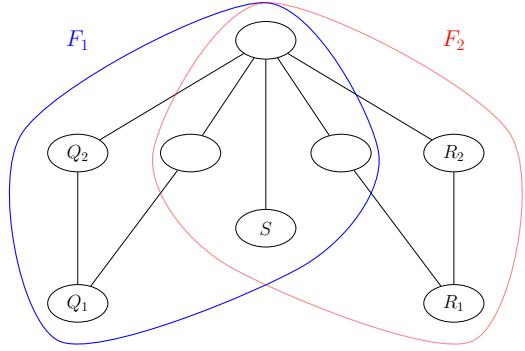

\begin{proof}
By Corollary~\ref{cor:O(P)_bijection}, the square faces of $\mathcal{O}(P)$ are in bijection with the set
\begin{equation*}
    S_{\mathcal{O}(P)} =
    \left\{ \{ F_1,F_2 \} \;\middle|\; 
        \begin{array}{l}
            F_1 \text{ and } F_2 \text{ filters such that} \\
            F_1 \setminus F_2 \text{ and } F_2 \setminus F_1 \text{ are connected}
        \end{array}
     \right\}.
\end{equation*}
By Corollary~\ref{cor:C(P)_bijection} the square faces of $\mathcal{C}(P)$ are in bijection with the set
\begin{equation*}
    S_{\mathcal{C}(P)} =
    \left\{ ( \{ Q,R \},S)  \;\middle|\;
    \begin{array}{l}
        Q = Q_1 \sqcup Q_2 \text{ connected, $Q_i$ antichains }  \\
        R = R_1 \sqcup R_2 \text{ connected, $R_i$ antichains } \\
        S \text{ antichain } \\
        Q \parallel S, \ Q \parallel R, \ R \parallel S
    \end{array}
    \right\}.
\end{equation*}
We aim to find maps $S_{\mathcal{O}(P)} \to S_{\mathcal{C}(P)}$ and $S_{\mathcal{C}(P)} \to S_{\mathcal{O}(P)}$ that are inverses of each other.

First we define the map
\begin{align*}
    \Phi \colon S_{\mathcal{O}(P)} &\longrightarrow S_{\mathcal{C}(P)} \\
    \{F_1, F_2\}
    &\longmapsto (\{ Q,R \},S),
\end{align*}
where
\begin{align*}
    Q &\coloneqq \min(F_1 \setminus F_2) \cup \max(F_1 \setminus F_2) , \\
    R &\coloneqq \min(F_2 \setminus F_1) \cup \max(F_2 \setminus F_1),  \\
    S &\coloneqq \{ p \in \min(F_1 \cap F_2 ) \colon p \parallel Q \text{ and } p \parallel R \}.
\end{align*}
See Figure \ref{fig:schematic} for a schematic picture of $\Phi$.
Clearly $S$ is an antichain parallel to $Q$ and $R$.
By Lemma~\ref{lm:extremeconn}, both $Q$ and $R$ are connected.
Here both $Q$ and $R$ can be written as a disjoint union of antichains (if, say, $Q = \{ p \}$ then we can write $Q = \emptyset \sqcup \{ p \}$).
Since $F_1$ and $F_2$ are filters we have $(F_1 \setminus F_2) \parallel (F_2 \setminus F_1)$ and thus $Q \parallel R$.
Therefore $\Phi$ is well-defined.

We then define the map
\begin{align*}
    \Psi \colon S_{\mathcal{C}(P)} &\longrightarrow S_{\mathcal{O}(P)} \\
    (\{ Q,R \},S) &\longmapsto 
    \left\{
    \begin{array}{l}
        F_1 \coloneqq \filter{Q} \cup \filter{S} \cup ( \filter{R} \cap \ideal{R}^\mathsf{c} ), \\
        F_2 \coloneqq \filter{R} \cup \filter{S} \cup (\filter{Q} \cap \ideal{Q}^\mathsf{c})
    \end{array}
    \right\},
\end{align*}
where for an ideal $J \subseteq P$ we write $J^\mathsf{c} \coloneqq P \setminus J$, which is a filter.
Figure~\ref{fig:schematic} gives also a schematic picture of the map $\Psi$.
Clearly $F_1$ and $F_2$ are filters.
Here
\begin{align*}
    F_1\setminus F_2 &= \filter{Q} \setminus (\filter{R} \cup \filter{S} \cup (\filter{Q} \cap \ideal{Q}^\mathsf{c})) \\
    &\cup \underbrace{\filter{S} \setminus  (\filter{R} \cup \filter{S} \cup (\filter{Q} \cap \ideal{Q}^\mathsf{c}))}_{= \emptyset}  \\
    &\cup \underbrace{(\filter{R} \cap \ideal{R}^\mathsf{c}) \setminus (\filter{R} \cup \filter{S} \cup (\filter{Q} \cap \ideal{Q}^\mathsf{c}))}_{= \emptyset} \\
    &= \filter{Q} \setminus (\filter{R} \cup \filter{S} \cup (\filter{Q} \cap \ideal{Q}^\mathsf{c})) \\
    &\stackrel{(\ast)}{=} \filter{Q} \cap \ideal{Q},
\end{align*}
where we prove $(\ast)$ as follows:
If $x \in \filter{Q} \setminus (\filter{R} \cup \filter{S} \cup (\filter{Q} \cap \ideal{Q}^\mathsf{c}))$ then 
\begin{equation*}
    x \in \filter{Q} \setminus (\filter{Q} \cap \ideal{Q}^\mathsf{c}) = \filter{Q} \cap \ideal{Q}.
\end{equation*}
If $x \in \filter{Q} \cap \ideal{Q}$ then we can't have $x \in \filter{R} \cup \filter{S}$ since $Q \parallel R$ and $Q \parallel S$.
Therefore $x \in \filter{Q} \setminus (\filter{R} \cup \filter{S} \cup (\filter{Q} \cap \ideal{Q}^\mathsf{c}))$.
Thus we have $(\ast)$.
By Lemma \ref{lemma:connected}, $F_1 \setminus F_2 = \filter{Q} \cap \ideal{Q}$ is connected.
By symmetry, we have $F_2 \setminus F_1 = \filter{R} \cap \ideal{R}$ which is also connected by Lemma \ref{lemma:connected}.
Thus $\Psi$ is well-defined.

We now turn to proving that $\Phi$ and $\Psi$ are inverses of each other

\medbreak
\textbf{$\Phi \circ \Psi = \id$:} 
Let $(\{ Q,R\},S) \in S_{\mathcal{C}(P)}$.
Let
\begin{equation*}
    \Psi(\{ Q,R \} ,S ) =
    \left\{
    \begin{array}{l}
        F_1 \coloneqq \filter{Q} \cup \filter{S} \cup ( \filter{R} \cap \ideal{R}^\mathsf{c} ), \\
        F_2 \coloneqq \filter{R} \cup \filter{S} \cup (\filter{Q} \cap \ideal{Q}^\mathsf{c})
    \end{array}
    \right\}.
\end{equation*}
We need to show that $(\{ Q,R\}, S ) = \Phi(\{ F_1,F_2 \})$, that is,
\begin{equation}
    \{ Q,R\} =
    \left\{ 
        \begin{array}{l}
            \min(F_1 \setminus F_2) \cup \max(F_1 \setminus F_2), \\
            \min(F_2 \setminus F_1) \cup \max(F_2 \setminus F_1) \\
        \end{array}
     \right\} \label{eq:QRS}
\end{equation}
and
\begin{equation}
    S = \{ p \in \min(F_1 \cap F_2) \colon p \parallel \min(F_i \setminus F_j) \cup \max(F_i \setminus F_j), \ i,j= 1,2 \}. \label{eq:S_equals}
\end{equation}
We saw above that $F_1 \setminus F_2 = \filter{Q} \cap \ideal{Q}$ and $F_2 \setminus F_1 = \filter{R} \cap \ideal{R}$.
Since $Q$ and $R$ have height $\le 2$, we have $\min(F_1 \setminus F_2) \cup \max(F_1 \setminus F_2) = Q$ and $\min(F_2 \setminus F_1) \cup \max(F_2 \setminus F_1) = R$ by Lemma~\ref{lemma:min_cup_max}.
It remains to show (\ref{eq:S_equals}).
Let $x \in S$.
Now $x \in F_1 \cap F_2$.
If we had $y < x$ for some $y \in F_1 \cap F_2$ then this would be a contradiction since $Q \parallel S$, $R \parallel S$ and $S$ is an antichain.
Therefore $x \in \min(F_1 \cap F_2)$.
Since $Q \parallel S$ and $R \parallel S$, $x$ cannot be comparable to any element in $F_1 \setminus F_2 = \filter{Q} \cap \ideal{Q}$ or to any element in $F_2 \setminus F_1 = \filter{R} \cap \ideal{R}$.
In particular, $x$ cannot be comparable to any element in $\max(F_i \setminus F_j)$ or $\min(F_i \setminus F_j)$.
Then suppose $x$ belongs to the set on the right-hand side in (\ref{eq:S_equals}).
Now $x \parallel Q$ and $x \parallel R$.
Since $x \in F_1 \cap F_2$, we need to have $x \in \filter{S}$.
Since $x$ is a minimal element in $F_1 \cap F_2$, we need to have $x \in S$.
This finishes our proof of $\Phi \circ \Psi = \id$.

\medbreak
\textbf{$\Psi \circ \Phi = \id$:}
Let $\{ F_1,F_2 \} \in S_{\mathcal{O}(P)}$.
Let $\Phi(\{ F_1,F_2 \}) = (\{ Q,R \},S)$
where
\begin{align*}
    Q &\coloneqq \min(F_1 \setminus F_2) \cup \max(F_1 \setminus F_2), \\
    R &\coloneqq \min(F_2 \setminus F_1) \cup \max(F_2 \setminus F_1), \\
    S &\coloneqq \{ p \in \min(F_1 \cap F_2) \colon p \parallel Q  \text{ and } p \parallel R\}.
\end{align*}
We aim to show that $\{ F_1,F_2 \} = \Psi(\{ Q,R,S \})$, that is,
\begin{equation}
    \{ F_1,F_2 \} = 
    \left\{
    \begin{array}{l}
         \filter{Q} \cup \filter{S} \cup ( \filter{R} \cap \ideal{R}^\mathsf{c} ), \\
         \filter{R} \cup \filter{S} \cup (\filter{Q} \cap \ideal{Q}^\mathsf{c})
    \end{array}
    \right\}. \label{eq:F1_F2}
\end{equation}
We claim that
\begin{equation}
    \filter{Q}  \cup \filter{S} \cup (\filter{R} \cap \ideal{R}^\mathsf{c}) = F_1. \label{eq:F1}
\end{equation}
First we prove $\subseteq$ in (\ref{eq:F1}).
Since $F_1$ is a filter and $S,Q \subseteq F_1$ we have $\filter{Q} \subseteq F_1$ and $\filter{S} \subseteq F_1$.
Let $x \in \filter{R} \cap \ideal{R}^\mathsf{c}$.
Now $x \ge r$ for some $r \in R \subseteq F_2$.
Thus $x \in F_2$.
If we had $x \not\in F_1$ then $x \in F_2 \setminus F_1$ and thus $x$ would be below some element from $\max(F_2 \setminus F_1)$ and hence $x \in \ideal{R}$, which is a contradiction.
Hence we must have $x \in F_1$.
We have now shown $\subseteq$ in (\ref{eq:F1}).
For the direction $\supseteq$ let $x \in F_1$.
If $x \not\in F_2$ then $x \in F_1 \setminus F_2$ and thus $x$ is above some element from $\min(F_1 \setminus F_2)$, which implies $x \in \filter{Q}$.
Thus suppose $x \in F_1 \cap F_2$.
Now $x \ge m$ for some $m \in \min(F_1 \cap F_2)$.
If $m$ is not comparable to any element in $Q$ or $R$, then $m \in S$ and thus $x \in \filter{S}$.
Thus suppose $m$ is comparable to some element in $Q \cup R$.
Suppose first that $m \perp q$ for some $q \in Q$. 
Since $Q \subseteq F_1 \setminus F_2$ and $m \in F_1 \cap F_2$, $m$ cannot be below any element from $Q$.
Therefore $m \ge q$ which implies $x \in \filter{Q}$.
Then suppose $m \perp r$ for some $r \in R$.
Similarly, we must have $m \ge r$ and thus $x \in \filter{R}$.
Since $R \subseteq F_2 \setminus F_1$ and $x \in F_1$, $x$ cannot be below any element from $R$.
Hence $x \in \ideal{R}^\mathsf{c}$ and therefore $x \in \filter{R} \cap \ideal{R}^\mathsf{c}$.
We have now shown (\ref{eq:F1}).
By symmetry, we have 
\begin{equation*}
    \filter{R}  \cup \filter{S} \cup (\filter{Q} \cap \ideal{Q}^\mathsf{c}) = F_2.
\end{equation*}
Hence we have (\ref{eq:F1_F2}).
We are now done showing $\Psi \circ \Phi = \id$, which finishes the proof of the theorem.
\end{proof}

\section{Conclusions}

In this paper, we have studied the structure of two-dimensional faces of the order polytope $\Oc(P)$ and the chain polytope $\Cc(P)$. By Corollary~\ref{cor:triangle_or_square}, we have $f_2(\Pc)=f_2^\Delta(\Pc)+f_2^\Box(\Pc)$ when $\Pc$ is the order polytope or the chain polytope of the poset $P$. Combining Theorem~\ref{thm:triangle_ineq} and Theorem~\ref{thm:square_eq}, we thus deduce the case $i=2$ for Conjecture~\ref{HL-conj}:

\begin{theorem}\label{main}
Let $P$ be a poset. Then 
\begin{equation*}
    f_2(\Oc(P))=f_2^\Delta(\Oc(P))+f_2^\Box(\Oc(P))\leq f_2^\Delta(\Cc(P))+f_2^\Box(\Cc(P)) =f_2(\Cc(P)),
\end{equation*}
with equality holding if and only if $P$ is $\mathbf{X}$-free.
\end{theorem}
We immediately get the following corollary, which adds to the list of equivalent conditions in~\cite[Corollary 3.2]{HLSS}.
\begin{corollary}
    Let $P$ be a poset. The following are equivalent.
    \begin{itemize}
      \item $\Oc(P)$ and $\Cc(P)$ are unimodularly equivalent.
      \item $f_2(\mathcal{O}(P)) = f_2(\mathcal{C}(P))$.
      \item $f_2^\Delta(\mathcal{O}(P)) = f_2^\Delta(\mathcal{C}(P))$.
    \end{itemize}
\end{corollary}

\subsection*{Acknowledgements}
This work was partially carried out while the third author was visiting Aalto University, supported by the Visiting Researcher Program of the Aalto Science Institute. 
We thank Aalto University and the Aalto Science Institute for their hospitality. 
We also thank the anonymous referee for carefully reading the paper and providing valuable comments and suggestions.

\end{document}